\documentclass{amsart}
\usepackage{amsfonts}
\usepackage{amssymb}
\usepackage{amsthm}
\usepackage{amsmath}
\usepackage[dvips]{graphics,color}
\usepackage{verbatim}

\newtheorem{thm}{Theorem}[section]
\newtheorem{defi}{Definition}[section]
\newtheorem{corol}{Corollary}[section]
\newtheorem{lemma}{Lemma}[section]
\newtheorem{prop}{Proposition}[section]

\newtheorem{claim}{Claim}

\newtheorem{problem}{Problem}[section]

\theoremstyle{definition}
\newtheorem{example}{Example}[section]
\newtheorem{remark}{Remark}[section]

\newcommand{\Proof}{\noindent {\bf Proof.~~}}
\newcommand{\eproof}{\hfill\rule{2.2mm}{3.0mm}}

\numberwithin{equation}{section}

\def\brho{\boldsymbol{\rho}}

\def\bN{\mathbb{N}}
\def\bQ{\mathbb{Q}}
\def\bR{\mathbb{R}}
\def\bZ{\mathbb{Z}}
\def\bZp{\mathbb{Z}^+}
\def\wdt{\widetilde}
\def\gF{\Phi}
\def\gY{\Psi}

\def\gt{\tau}

\def\ga{\alpha}
\def\gb{\beta}
\def\gd{\delta}
\def\gm{\mu}
\def\gn{\nu}
\def\gp{\pi}
\def\gs{\sigma}
\def\gl{\lambda}
\def\gL{\Lambda}
\def\gr{\rho}
\def\gS{\Sigma}
\def\gga{\gamma}
\def\gGa{\Gamma}
\def\sD{\mathcal{D}}
\def\sC{\mathcal{C}}

\def\sA{\mathcal{A}}

\def\sF{\mathcal{F}}
\def\sG{\mathcal{G}}
\def\sI{\mathcal{I}}
\def\sJ{\mathcal{J}}
\def\sM{\mathcal{M}}
\def\sP{\mathcal{P}}
\def\sS{\mathcal{S}}
\def\sT{\mathcal{T}}
\def\i{\mathbf{i}}
\def\j{\mathbf{j}}
\def\k{\mathbf{k}}
\def\u{\mathbf{u}}
\def\v{\mathbf{v}}

\def\sH{\mathcal{H}}
\def\sE{\mathcal{E}}

\def\diam{\mathrm{diam}\,}
\def\card{\textrm{card\,}}

\def\Tone{T^{(1)}}
\def\Done{D^{(1)}}
\def\Ttwo{T^{(2)}}
\def\Dtwo{D^{(2)}}
\def\Tthree{T^{(3)}}
\def\Dthree{D^{(3)}}
\def\Tfour{T^{(4)}}
\def\Dfour{D^{(4)}}

\def\sTOne{\sT^{(1)}}

\def\sTThree{\sT^{(3)}}

\def\empwd{{{\O}}}

\begin{document}           

\title{Lipschitz equivalence of self-similar sets with touching structures}  
\author{Huo-Jun Ruan \and Yang Wang \and Li-Feng Xi}

\begin{abstract}
   Lipschitz equivalence of self-similar sets is an important area in the study of fractal
   geometry. It is known that two dust-like self-similar sets with the same
   contraction ratios are always Lipschitz equivalent. However, when self-similar sets
   have touching structures the problem of Lipschitz equivalence becomes much more
   challenging and intriguing at the same time. So far the only known results only
   cover self-similar sets in $\bR$ with no more than 3 branches. In this study
   we establish results for the Lipschitz equivalence of self-similar sets with
   touching structures in $\bR$
   with arbitrarily many branches. Key to our study is the introduction
   of a geometric condition for self-similar sets called {\em substitutable}.
\end{abstract}
\maketitle

\section{Introduction}

\subsection{\emph{Motivation}}

A fundamental concept in fractal geometry is dimension.
It is often used to differentiate fractal sets, and when two sets
have different dimensions (Hausdorff dimensions, box dimensions or other dimensions)
we often consider them to be not alike. However two compact sets, even with the
same dimension, may in fact be quite different in many ways. Thus it is natural to
seek a suitable quality that would allow us to tell whether two
fractal sets are similar. And generally, Lipschitz
equivalence is thought to be such a suitable quality.

It has been pointed out in \cite{FaMa92} that while topology may be regarded as the study
of equivalence classes of sets under homeomorphism, fractal
geometry is sometimes thought of as the study of equivalence classes
under bi-Lipschitz mappings. The more restrictive maps such as isometry
tend to lead to poor and rather boring equivalent classes, while the far less restrictive
maps such as general continuous maps take us completely out of geometry
into the realm of pure topology (see~\cite{Gromo07}). Bi-Lipschitz maps offer
a good balance, which lead to categories that are interesting and intriguing
both geometrically and algebraically.

There are many works done in the filed of Lipschitz equivalence of two fractal sets. Some earlier fundamental results are obtained by Cooper and Pignataro \cite{CP88}, David and Semmes \cite{DaSe97}, and Falconer and Marsh \cite{FaMa89, FaMa92}. Recently, based on these works and motivated by Problem~11.16 in \cite{DaSe97}, Rao, Ruan, Wang, Xi, Xiong and their collaborators obtained a series of results, see e.g. \cite{RRW}-\cite{RRY08}, \cite{Xi10}-\cite{XiXi12}. There are also some other related works. Xi \cite{Xi04} discussed the nearly Lipscchitz equivalence of self-conformal sets. Mattila and Saaranen \cite{MaSa09} studied the Lipschitz equivalence of Ahlfors-David regular sets. Deng, Wen, Xiong and Xi \cite{DWXX11} and Llorente and Mattila \cite{LloMat10} discussed the bi-Lipschitz embedding of fractal sets.

Let $E$ and $F$ be two compact subsets  of $\bR^d$. A bijection
$f:\, E\to F$ is said to be \emph{bi-Lipschitz} if there exist two
positive constants $c$ and $c^\prime$ such that
\begin{equation}\label{eq:f-def-temp}
  c|x-y|\leq |f(x)-f(y)|\leq c'|x-y|, \quad \forall x,y\in E.
\end{equation}
$E$ and $F$ are said to be \emph{Lipschitz equivalent}, denoted by
$E\sim F$, if there exists a bi-Lipschitz map $f$ from $E$ to $F$.

We recall some basic notations in fractal geometry. Given a
family of similitude $\gF_i(x)$, $i=1,\ldots,n$, on $\bR^d$, where
each $\gF_i$ has contraction ratio $\rho_i$ with $\rho_i<1$, there exists a
unique nonempty compact subset $E$ of $\bR^d$ such that
$\bigcup_{i=1}^n\gF_i(E)=E$, see \cite{Hut81}. The set of maps
$\{\gF_i(x), i=1,\ldots,n\}$ is called an \emph{iterated
function system} (IFS) and $E$ is called the
\emph{attractor}, or the \emph{invariant set}, of the IFS. We
also call $E$ a \emph{self-similar set} since every $\gF_i$ is a
similitude. If $\gF_i(E)\cap \gF_j(E)=\emptyset$ for any distinct
$i$ and $j$, the IFS $\{\gF_i\}$ is then said to satisfy the \emph{strong
separation condition (SSC)}, and $E$ is said to be \emph{dust-like}.

Given $\rho_1,\ldots,\rho_n\in (0,1)$ with $\sum_{i=1}^n
\rho_i^d<1$, we call $\brho=(\rho_1,\ldots,\rho_n)$ a (separable)
{\em contraction vector} (in $\mathbb{R}^d$). We denote by $\sD(\brho)$
the family of all dust-like self-similar sets with contraction
vector $\brho$ (here the ambient dimension $d$ is implicitly fixed).
The following property is well known, see e.g.
\cite{RRX06}.

\begin{prop}\label{prop1}
  Any two sets in $\sD(\brho)$ are Lipschitz equivalent.
\end{prop}

There are examples where two different contraction vectors $\brho_1$ and $\brho_2$
lead to Lipschitz equivalent families $\sD(\brho_1)$ and $\sD(\brho_2)$. For example,
Rao, Ruan and Wang \cite{RRW} has completely classified the Lipschitz equivalence of
dust-like families $\sD(\brho)$ where $\brho =(\rho_1, \rho_2)$, and one of the
results is that $\brho_1 =(\lambda, \lambda^5)$ and $\brho_2 =(\lambda^2, \lambda^3)$
lead to Lipschitz equivalence families whenever the resulting self-similar
sets are dust-like. The paper \cite{RRW} and some earlier studies such as \cite{CP88,FaMa92}
have explored the impact of algebraic properties of the contraction vectors on
Lipschitz equivalence, yielding a number of intriguing results showing the links.

Nevertheless one should not overlook the importance of geometric properties of
the underlying IFSs has on Lipschitz equivalence of self-similar sets. Relating to this
point is an interesting problem proposed by David and Semmes \cite{DaSe97}
(Problem~11.16).
\begin{problem}\label{prob1}
  Let $S_i(x):=x/5+(i-1)/5$ be a contractive map from
  $[0,1]$ to $[0,1]$ where $i\in \{1,\cdots,5\}$. Let $M$ and
  $M^\prime$ be the attractor of the IFS $\{S_1,S_3,S_5\}$ and
  the IFS $\{S_1,S_4,S_5\}$, respectively. Are $M$ and $M^\prime$ Lipschitz
  equivalent?
\end{problem}

\begin{figure}[htbp]
    \center
  \scalebox{0.5}{\includegraphics{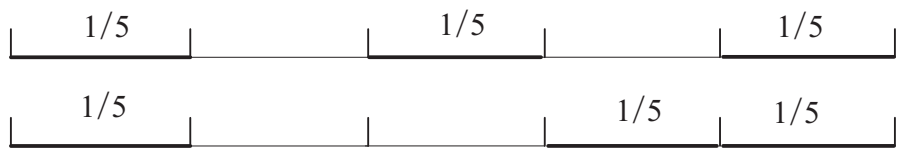}}
  \caption{Initial construction of $M$ and $M^\prime$}
  \label{figure: DS-set}
\end{figure}

We call $M$ the $\{1,3,5\}$-set and $M^\prime$ the $\{1,4,5\}$-set.
Clearly, $M$ is dust-like and $M^\prime$ has certain touching
structure, see Figure~\ref{figure: DS-set}. In this problem, the contraction
rations are all identical so the difference lies entirely in the geometry of
the two IFSs. David and Semmes
conjectured that $M\not\sim M^\prime$. However, by examining graph-directed
structures of the attractors and introducing techniques to study
Lipschitz equivalence on these structures, Rao, Ruan and Xi
\cite{RRX06} proved that in fact $M\sim M^\prime$.

Notice that all contractive maps in above problem have same contraction ratio $1/5$. Some similar works has been done in higher dimensional case, e.g. \cite{LauLuo12,Roi10,WZD12,XiXi10}.

\begin{figure}[htbp]
    \center
  \scalebox{0.5}{\includegraphics{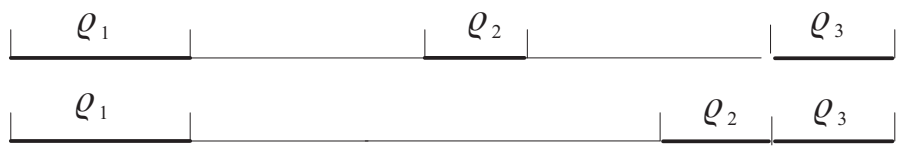}}
  \caption{Initial construction of $M_{\brho}$ and
  $M^\prime_{\brho}$}
  \label{figure: general-DS-set}
\end{figure}

A follow up study in Xi and Ruan \cite{XiRu07} exploits the interplay of algebraic
properties of contraction vectors and geometric properties of IFSs. It
considers the following generalization of the $\{1,3,5\}-\{1,4,5\}$ problem:

\begin{problem}\label{prob2}
  Let $\brho=(\gr_1,\gr_2,\gr_3)$ be a contraction vector
  (in $\mathbb{R}$). Let $\gF_i(x)=\gr_i x+d_i$, $i\in \{1,2,3\}$, where
  $d_1=0, d_3=1-\gr_3$ and $\gr_1<d_2<1-\gr_2-\rho_3$ (e.g.
  $d_2=\gr_1+(1-\gr_1-\gr_2-\gr_3)/2$).
  Let $\gY_1=\gF_1$, $\gY_3=\gF_3$ and $\gY_2(x)=\gr_2 x+t_2$ with
  $t_2=1-\gr_2-\gr_3$.
  Let $M_{\brho}$ and
  $M^\prime_{\brho}$ be the attractor
  of $\{\gF_1,\gF_2,\gF_3\}$  and
  $\{\gY_1,\gY_2,\gY_3\}$, respectively. See
  Figure~\ref{figure: general-DS-set} for their initial configuration. Are $M_{\brho}$
   and $M^\prime_{\brho}$ Lipschitz
  equivalent?
\end{problem}
Somewhat surprisingly, the Lipschitz equivalence of the two sets are completely
determined by the algebraic property of $\gr_1$ and $\gr_3$ and independent
of  $\gr_2$. It is shown in \cite{XiRu07}
that $M_{\brho}\sim
M^\prime_{\brho}$ if and only if $\log
\gr_1/\log\gr_3\in \bQ$.

The above example is nevertheless a very special case. It is natural to exploit such
algebraic and geometric connections further in more
general settings, which is the aim of this paper. Given the complexity of even to
establish the result for Problem \ref{prob2}, this may appears to be a very daunting
task. Fortunately, by introducing a new geometric notion called {\em substitutable}
we are able to prove a number of results in this direction.

Throughout this paper we assume that $\brho=(\rho_1,\ldots,\rho_n)$ is a
contraction vector (in $\bR$) with $n\geq 3$.  Let $D\in \sD(\brho)$. By Proposition~\ref{prop1},
we may assume without loss of generality that $D$ is the attractor
of the IFS $\{\gF_i(x)=\rho_i x+d_i\}_{i=1}^n$, where
$\gF_1([0,1]),\ldots,\gF_n([0,1])$ are equally spaced closed subintervals
of $[0,1]$ arranged from left to right,
normalized so that the left endpoint of $\gF_1([0,1])$ is $0$ and the right
end of $\gF_n([0,1])$ is $1$.

We are interested in the Lipschitz
equivalence of $D$ with the attractor $T$ of another IFS
$\{\gY_i(x)=\rho_i x+t_i\}_{i=1}^n$ having the same contraction vector $\brho$
but with translations $\{t_i\}$ that may result in some of the subintervals
$\gY_1([0,1]), \dots, \gY_n([0,1])$ touching one another (but no overlapping).
More precisely, the IFS $\{\gY_i(x)=\rho_i x+t_i\}_{i=1}^n$ satisfies
the following three properties:
\begin{itemize}
\item[\rm (1)]~The subintervals
$\gY_1([0,1]),\ldots,\gY_n([0,1])$ are spaced from left to right
without overlapping, i.e. their interiors do not intersect.
\item[\rm (2)]~The left endpoint of $\gY_1[0,1]$ is $0$ and the right endpoint
of $\gY_n[0,1]$ is $1$.
\item[\rm (3)]~There exists at least one $i\in\{1,2,\ldots,n-1\}$, such that the intervals
$\gY_i([0,1])$ and $\gY_{i+1}([0,1])$ are touching, i.e.
$\gY_i(1)=\gY_{i+1}(0)$.
\end{itemize}

\begin{figure}[htbp]
    \center
  \scalebox{0.5}{\includegraphics{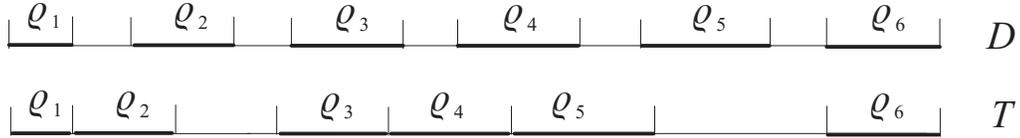}}
  \caption{Initial construction of $D$ and $T$, where $n=6$}
  \label{figure: n-branch-set}
\end{figure}

Denote by $T$ the attractor of the IFS $\{\gY_i\}_{i=1}^n$.
Figure~\ref{figure: n-branch-set} gives an example of $\{\gF_i\}$
and $\{\gY_i\}$, respectively. In this paper we present
necessary conditions and sufficient conditions for $D\sim T$.

\subsection{\emph{Notations and Examples}}

First some commonly used basic notations.
Denote $\Sigma_n:=\{1, 2, \dots,n\}$ and
$\Sigma_n^*:= \bigcup_{m\geq 1}\Sigma_n^m = \bigcup_{m\geq 1}\{1,2,\dots,n\}^m$.
We shall call any $i\in\Sigma_n$ a {\em letter} and
$\i=i_1\cdots i_m\in\Sigma_n^*$ a {\em word} of length
$|\i|:=m$. $i_1$ and $i_m$ is called the {\em first letter} and {\em last letter} of $\i$, respectively. We define
$\gr_\i=\gr_{i_1}\cdots\gr_{i_m}$,
$\gY_\i=\gY_{i_1}\circ\cdots\circ\gY_{i_m}$ and $T_\i=\gY_\i(T)$.
$T_\i$ is called a {\em cylinder} of the IFS $\{\gY_i\}$ for $\i$.
Similarly we define
$\gF_\i=\gF_{i_1}\circ\cdots\circ\gF_{i_m}$ and the cylinder
$D_\i=\gF_\i(D)$.

Specific to this study we introduce also other notations. A letter $i\in\Sigma_n$
is a {\em (left) touching letter} if $\gY_i([0,1])$ and $\gY_{i+1}([0,1])$
are touching, i.e. $\gY_i(1) =\gY_{i+1}(0)$. We use
$\Sigma_T\subset\Sigma_n$ to
denote the set of all (left) touching letters. Note that one may view
$\Sigma_T+1$ to be the set of all right touching letters. For simplicity
we shall drop the word ``left'' for $\Sigma_T$.
Let $\ga$ and $\gb$ be the number of successive touching intervals
among $\gY_1([0,1]),\ldots,\gY_n([0,1])$ at the beginning and at the end, respectively.
In other words, $\bigcup_{i=1}^\ga \gY_i[0,1]$ and $\bigcup_{i=n-\gb+1}^n
\gY_i[0,1]$ are intervals, while $\gY_\ga(1)\not=\gY_{\ga+1}(0)$ and
$\gY_{n-\gb}(1)\not= \gY_{n-\gb+1}(0)$.

Given a cylinder $T_\i$ and a nonnegative integer $k$, we can
define respectively the {\em level $(k+1)$ left touching patch} and the
{\em level $(k+1)$ right touching patch}
of $T_\i$ to be
\begin{equation}\label{eq:Lk-Rk-def}
  L_k(T_\i)=\bigcup_{j=1}^{\ga} T_{\i[1]^k j}, \quad
  R_k(T_\i)=\bigcup_{j=n-\gb+1}^n T_{\i [n]^k j},
\end{equation}
where $[\ell]^k$ is defined to be the word $\underbrace{\ell\cdots\ell}_k$
for any $\ell\in \{1,\ldots,n\}$, with $\i[1]^k j$ be the
concatenation of $\i$, $[1]^k$ and the letter $j$ (similarly for
$\i [n]^k j$). We remark that
$L_0(T_\i)=\bigcup_{j=1}^{\ga} T_{\i j}$ and $R_0(T_\i)=\bigcup_{j=n-\gb+1}^n T_{\i j}$.

Now comes the main notation we introduce for this paper.
A letter $i\in \Sigma_T$ is called \emph{left substitutable} if there exist $\j\in
\Sigma_n^*$ and $k,k'\in \bN$, such that $\diam L_{k}(T_{i+1})=\diam
L_{k'}(T_{i\j})$ and the last letter of $\j$
does not belong to $\{1\}\cup(\gS_T+1)$. Geometrically it simply means that
certain left touching patch of the cylinder $T_{i+1}$ has the same diameter as
that of some left touching patch
of a cylinder $T_{i\j}$, and as a result we can substitute one of the left
touching patches by the
other without disturbing the other neighboring structures in $T$ because they
have the same diameter. The actual substitution is performed in the proof of our
main theorem. Similarly, $i\in \Sigma_T$
is called \emph{right substitutable} if there exist $\j\in
\Sigma_n^*$ and $k,k'\in \bN$, such that $\diam R_{k}(T_{i})=\diam
R_{k'}(T_{(i+1)\j})$ and the last letter of $\j$
does not belong to $\{n\}\cup\gS_T$. We say that $i\in\Sigma_T$ is
\emph{substitutable} if it is left substitutable or right
substitutable.

\begin{remark}
  Both left and right substitutable properties can be characterized algebraically as well.
  By definition, it is easy to check that
  $\diam L_{k}(T_{i+1})=\diam L_{k'}(T_{i\j})$ is
  equivalent to
  \begin{equation}\label{eq:left-disp}
    \gr_{i+1}\gr_1^k=\gr_i\gr_1^{k^\prime}\gr_\j,
  \end{equation}
  while $\diam R_{k}(T_{i})=\diam R_{k'}(T_{(i+1)\j})$ is equivalent
  to
  \begin{equation}\label{eq:right-disp}
    \gr_{i}\gr_n^k=\gr_{i+1}\gr_n^{k^\prime}\gr_\j.
  \end{equation}
\end{remark}

\begin{figure}[htbp]
    \center
  \scalebox{0.5}{\includegraphics{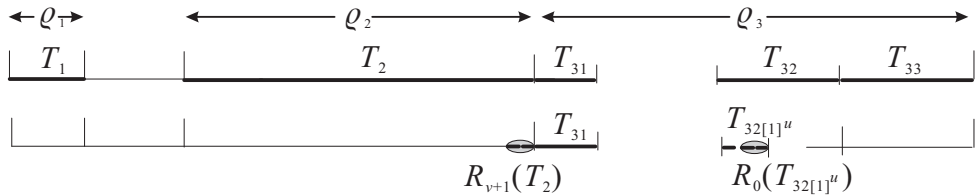}}
  \caption{The unique touching letter $2$ is right substitutable in Example~\ref{exam2}}
  \label{figure: general-DS-set-substitutable}
\end{figure}

\begin{example}\label{exam2}
  Let $\gY_1, \gY_2,\gY_3$ be defined as in Problem~\ref{prob2} and let
  $T$ be its attractor. Clearly $\gS_T=\{2\}$, $\ga=1$ and $\gb=2$. Assume that $\log
  \gr_1/\log\gr_3\in\bQ$, i.e. there exist $u,v\in\bZ^+$ such that
  $\gr_1^u=\gr_3^v$. Pick $k=v+1$, $k^\prime=0$ and $\j=2[1]^u$. It
  is easy to check that (\ref{eq:right-disp}) holds for $i=2$ and the
  last letter of $\j$ is
  $1\not\in\{3\}\cup\gS_T$. Thus the touching letter $2$ is right substitutable.  See
  Figure~\ref{figure: general-DS-set-substitutable} for a graphical
  illustration.
\end{example}

\subsection{\emph{Statement of Results}}

We establish several results in this paper. First we
prove the following necessary condition for $D\sim T$, regardless of the geometric configuration of the IFS $\{\gY_i\}$:
\begin{thm}\label{thm:necessary-cond}
  Assume that $D\sim T$. Then $\log \gr_1/\log\gr_n\in\bQ$.
\end{thm}

As a result we shall always assume in this paper that $\log
\gr_1/\log\gr_n\in\bQ$. For the case of $n=3$ branches it was shown
in \cite{XiRu07} that the condition $\log \gr_1/\log\gr_3\in\bQ$ is
also sufficient for the Lipschitz equivalence of $D$ and $T$. So
naturally one may ask whether this condition is sufficient in
general. The following theorem shows that this is false, even for
the 4-branch case.

\begin{thm}\label{thm:branch-4}
  Let $n=4$, $\gr_1=\gr_4$, and $\gS_T=\{2\}$. Assume that $D\sim T$.
  Let $s$ be the common Hausdorff dimension of $D$ and $T$ and $\gm_i=\gr_i^s$ for $1\leq i \leq 4$.
  Then $\gm_2$ and $\gm_3$ must be algebraically dependent, namely there exists a
  rational nonzero polynomial $P(x,y)$ such that $P(\gm_2,\gm_3)=0$.
\end{thm}

Later in the paper we shall see that if $\log \gr_i/\log\gr_j \in \bQ$ for all
$i,j\in \{1,\ldots,n\}$ then $D\sim T$. To go deeper we must take into
account the geometric information of the IFS $\{\gY_i\}$. The main theorem of the
paper is:

\begin{thm}\label{thm:sufficient-cond}
  Assume that $\log \gr_1/\log \gr_n\in\bQ$. Then, $D\sim T$ if
  every touching letter for $T$ is substitutable.
\end{thm}

As indicated in Example~\ref{exam2}, the IFS in question in Problems \ref{prob2} has a single touching letter $\Sigma_T =\{2\}$ and
the letter is substitutable. Thus the Lipschitz equivalences in both problems follow directly from Theorem~\ref{thm:sufficient-cond}.

\begin{corol}[\cite{XiRu07}]
  Let $\brho = (\rho_1, \rho_2, \rho_3)$ and $M_{\brho}$
   and $M^\prime_{\brho}$ be sets defined in
   Problem~\ref{prob2}.  Then $M_{\brho}
   \sim M^\prime_{\brho}$ if and only if
   $\log\gr_1/\log\gr_3\in\bQ$.
\end{corol}

\begin{corol}[\cite{RRX06}]
  Let $M$ and $M^\prime$ be sets defined in Problem~\ref{prob1}.
  Then $M\sim M^\prime$. 
\end{corol}


Theorem~\ref{thm:sufficient-cond} allows us to establish a
more general corollary. The argument used to show the substitutability in
Example \ref{exam2} is easily extended to prove the following corollary:

\begin{corol}  \label{corol-1.3}
  $D\sim T$ if one of the following conditions holds:
  \begin{enumerate}
    \item[\rm (1)] ~$\log \gr_i/\log \gr_j\in\bQ$ for
  all $i,j\in\{1,n,\ga\}\cup (\gS_T+1)$.

  \item[\rm (2)] ~$\log \gr_i/\log \gr_j\in\bQ$ for
  all $i,j\in\{1,n,n-\gb+1\}\cup \gS_T$.
  \end{enumerate}
\end{corol}
\Proof~~Without loss of generality, we only prove (1). Given a
touching letter $i$. We will show that $i$ is right substitutable.
Since $\log \gr_{\ga}/\log \gr_n, \log \gr_{i+1}/\log \gr_n\in\bQ$,
there exist $u,v,w\in\bZ^+$ such that
$\gr_{\ga}^u=\gr_n^v=\gr_{i+1}^w$. Pick $k=2v$, $k^\prime=0$ and
$\j=i[i+1]^{w-1}[\ga]^{u}$. It is easy to check that
(\ref{eq:right-disp}) holds. Notice that
$\ga\not\in\{n\}\cup\gS_T$. It follows that
 $i$ is right substitutable.\eproof

\medskip

The following result, which we wish to state as a theorem because of the
simplicity of its statement, is a direct corollary of Corollary \ref{corol-1.3}.

\begin{thm}\label{thm:suff-cond-log}
  Assume that $\log \gr_i/\log\gr_j \in \bQ$ for all
  $i,j\in \{1,\ldots,n\}$. Then $D\sim T$.
\end{thm}

We remark that the above condition is clearly not a necessary condition, as we have seen
from the 3-branch case, for which the contraction ratio of the middle branch is
irrelevant. One difference between the dust-like case and the touching case is
that the order of the contraction ratios do matter, as Theorem~\ref{thm:necessary-cond}
indicates. However, the condition in Theorem \ref{thm:suff-cond-log} can be viewed as a
weak necessary condition in the sense that given a set of contraction ratios $\rho_1,
\dots, \rho_n$, if $\log \gr_i/\log\gr_j \not\in \bQ$ for some $i,j\in \{1,\ldots,n\}$ then
there exists a touching IFS whose contraction ratios are $\{\rho_i\}$ such that its
attractor $T$ is not Lipschitz equivalent to $D$. This is easily done by making the
contraction ratios of the left most and right most branches to be $\rho_i$ and $\rho_j$,
respectively.

%
%
%


The rest of the paper will be devoted to proving the stated results. In Section 2 we prove
Theorems~\ref{thm:necessary-cond} and \ref{thm:branch-4}, and in Section 3 we prove
Theorem \ref{thm:sufficient-cond}.

\section{Necessary condition for $D\sim T$}

\subsection{\emph{Bi-Lipschitz map related with a dust-like self-similar set}}

In this subsection, we will discuss the property of bi-Lipschitz map
$f: E\to F$, where $E$ is a nonempty compact subset of $\bR^d$ and
$F$ is a dust-like self-similar subset of $\bR^d$ with contraction
vector $(\rho_1,\ldots,\rho_n)$. We assume that
\begin{equation}\label{eq:f-def-temp}
  c|x-y|\leq |f(x)-f(y)|\leq c'|x-y|, \quad \forall x,y\in E,
\end{equation}
where $0<c\leq c'$.

For any nonempty subsets $A,B\subset \bR^d$, we define
$d(A,B)=\inf\{|x-y|:\, x\in A, y\in B\}$. The diameter of $A$ is
defined to be $\diam A:=\sup\{|x-y|:\, x,y\in A\}$. If $d(A,
B\setminus A)>0$, we say that $A$ is a \emph{$B$-separate set}. If
$d(A, B\setminus A)\geq \gl \cdot\diam A$ for some $\gl>0$, we say
that $A$ is a \emph{$(B,\gl)$-separate set}.

We now present a lemma which is similar to \cite[Lemma~3.2]{FaMa92}.

\begin{lemma}\label{lem:finite}
For any $\gl>0$, there exists an integer
  $n_0$ such that for any $(E,\gl)$-separate set $A\subset E$, there exist
$\k, \j_1, \dots, \j_p \in\Sigma_n^*$ such that $F_{\k\j_1}, \dots,
F_{\k\j_p}$ are disjoint and
\begin{equation} \label{3.1}
   f(A)=\bigcup_{r=1}^p F_{\k\j_r}\subset F_\k,
\end{equation}
where  each $|\j_r|= n_0$.
\end{lemma}

\Proof
  Given an $(E,\gl)$-separate set $A\subset E$. Let $F_\k$ be the
  smallest cylinder containing $f(A)$. Then, it is clear that there exists a positive constant $\gd$
  dependent only on $F$ such that $\diam F_\k\leq \gd \,\diam f(A)$.
  For a detailed proof, please see e.g. \cite[Lemma~3.1]{FaMa92}.
  Thus, by (\ref{eq:f-def-temp}), we have
  \begin{equation*}
    \diam F_\k\leq \gd \,\diam f(A)\leq \gd c'\diam A.
  \end{equation*}
  Let $n_0$ be the smallest integer satisfying $\overline{\rho}^{n_0}\gd
  c'\leq\frac{c\gl}{2}$, where
  $\overline{\rho}=\max\{\rho_1,\ldots,\rho_n\}$. Then, if
  $|\j'|=n_0$,
  \begin{equation}\label{eq:lem-finite-1}
    \diam F_{\k\j'} \leq \overline{\rho}^{n_0}\diam F_\k\leq  \overline{\rho}^{n_0}\gd
    c'\diam A\leq\frac{c\gl}{2}\diam A.
  \end{equation}

  Assume that $F_{\k\j'}\cap f(A)\not=\emptyset$, we will prove $F_{\k\j'}\subset
  f(A)$ by showing $d(z,f(E\setminus A))>0$ for any $z\in
  F_{\k\j'}$.

  Pick $x\in f^{-1}(F_{\k\j'}\cap f(A))$, we have
  \begin{equation}\label{eq:lem-finite-2}
    d(f(x),f(E\setminus A)) \geq d(f(A),f(E\setminus A)) \geq c\cdot
    d(A, E\setminus A) \geq c \gl\,\diam A.
  \end{equation}
  Thus, for any $z\in  F_{\k\j'}$,  using (\ref{eq:lem-finite-1}) and
  (\ref{eq:lem-finite-2}), we have
  \begin{equation*}
    d(z,f(E\setminus A)) \geq d(f(x),f(E\setminus A)) - \diam F_{\k\j'}
    \geq c \gl\, \diam A - \frac{c\gl}{2}\diam A>0.
  \end{equation*}
  This completes the proof of the lemma.
\eproof

\begin{remark}\label{rem:optimal-decomp}
By the proof, we can require that $F_\k$ is the smallest cylinder
containing $f(A)$. Under this restriction, $\k$ is uniquely
determined by $A$. Consequently, the set $\{\j_1,\ldots,\j_p\}$ are
also uniquely determined by $A$ and $n_0$.
\end{remark}

\subsection{\emph{Construction of $(T,\gl)$-separate sets}}

Let $\empwd$ be the empty word. We say that the length of $\empwd$ is $0$.
Define $\gY_\empwd=\gF_\empwd=\textrm{id}, \gr_\empwd=1$ and $I_\empwd=[0,1]$.
Let $\sI_0=\{I_\empwd\}$ and $\sI_m=\{I_\i: |\i|=m\}$ for all positive integers $m$.

By the definition of $T$, we know that there exists $i$ such that
$\gY_i(1)=\gY_{i+1}(0)$. We pick one such $i$ and denote it by
$i_0$. Without loss of generality, we assume that $\rho_1\geq
\rho_n$. For positive integer $k$, we define $\gt(k)$ to be the
unique positive integer satisfying
\begin{equation}\label{eq:gt-k-def}
  \rho_n^k \rho_1< \rho_1^{\gt(k)} \leq \rho_n^k.
\end{equation}
It is clear that $\gt(k)\geq k$ and is increasing with respect to
$k$. We define
\begin{equation*}
   C^k=R_k(T_{i_0}) \cup L_{\gt(k)}(T_{i_0+1}).
\end{equation*}
We remark that $C^1\supset C^2\supset\cdots$.

We shall adopt the notation $\asymp$ throughtout this paper. Let $A$ be a given index
set. Given two sequences of positive real numbers  $(a_i)_{i\in A}$ and $(b_i)_{i\in A}$
indexed by $A$, we denote $(a_i)\asymp (b_i)$ if there exist positive constants $c_1,
c_2$ independent of $i$ such that $c_1 a_i\leq b_i\leq c_2 a_i$ for all $i\in A$. For
convenience of statement in the proofs, we shall often write {\em $a_i\asymp b_i$ for
all $i\in A$}, or simply $a_i\asymp b_i$ if there is no confusion about the index set.

\begin{lemma}\label{lem:T-gl-separated}
  There exists $\gl>0$, such that $C^k$ is $(T,\gl)$-separate for
  all $k$.
\end{lemma}
\Proof
  Notice that $R_k(T_{i_0}) \cap L_{\gt(k)}(T_{i_0+1})$ is a singleton. By (\ref{eq:gt-k-def}),
  for all $k$ we have
  \begin{eqnarray}
    \diam C^k&=&\diam R_k(T_{i_0}) + \diam L_{\gt(k)}(T_{i_0+1}) \notag \\
    &=&\rho_{i_0} \rho_n^k \cdot \diam \Big(\bigcup_{j=n-\gb+1}^n
    T_j\Big) + \rho_{i_0+1} \rho_1^{\gt(k)} \cdot \diam
    \Big(\bigcup_{j=1}^\ga T_j\Big) \notag\\
    &\asymp& \rho_n^k. \label{eq:diam-Bk-asymp}
  \end{eqnarray}
On the other hand, it is clear that the distance of $C^k$ and $T\setminus
  C^k$ equals the minimum of the following two distances:
  $d(R_k(T_{i_0}), T_{i_0[n]^k(n-\gb)})$ and $d(L_{\gt(k)}(T_{i_0+1}),
  T_{(i_0+1)[1]^{\gt(k)}(\ga+1)})$. Since
  \begin{eqnarray*}
    d(R_k(T_{i_0}), T_{i_0[n]^k(n-\gb)}) = d(T_{i_0[n]^k(n-\gb+1)},
    T_{i_0[n]^k(n-\gb)}) = \rho_{i_0} \rho_n^k\cdot d(T_{n-\gb+1},
    T_{n-\gb}),
  \end{eqnarray*}
  and similarly,
  \begin{eqnarray*}
    d(L_{\gt(k)}(T_{i_0+1}),
  T_{(i_0+1)[1]^{\gt(k)}(\ga+1)})  = \rho_{i_0+1} \rho_1^{\gt(k)}\cdot d(T_\ga,
    T_{\ga+1}),
  \end{eqnarray*}
  we know from (\ref{eq:gt-k-def}) that $d(C^k, T\setminus C^k) \asymp
  \rho_n^k$ for all $k$. Combining this with (\ref{eq:diam-Bk-asymp}), we have
  $d(C^k, T\setminus C^k) \asymp \diam C^k$ for all $k$. Hence, there exists
  $\gl>0$ such that $d(C^k, T\setminus C^k)\geq \gl\cdot  \diam
  C^k$ for any $k$. This completes the proof.
\eproof

For all $\i\in\Sigma_n^*\cup \{\empwd\}$ and $k\in\bZp$, we define
\begin{equation*}
  C_{\i}^k=\gY_\i(C^k).
\end{equation*}
It is clear that $\diam C_{\i}^k = \rho_{\i} \cdot \diam C^k$ and
$d(C_{\i}^k, T\setminus C_{\i}^k) = \rho_{\i} \cdot d(C^k,
T\setminus C^k)$. Thus, $C_{\i}^k$ is $(T,\gl)$-separated, where
$\gl$ is defined as in Lemma~\ref{lem:T-gl-separated}. For any $k\in \bZp$, we define
$$\sC_k=\{C_\i^j:\, |\i|+j=k \textrm{ where $ \i\in \gS_n^*\cup\{\empwd\} $ and $j\in \bZp$}\}.$$

\begin{lemma}\label{lem:C-ik-disjoint}
  For any two distinct sets $A,B\in \sC_k$, we have $A\cap B=\emptyset$.
\end{lemma}
\Proof
  Suppose that $A=C_\i^k$ and $B=C_\j^\ell$.
  It is clear that $A\cap B =\emptyset$ if $k=\ell$. Thus, without loss of generality, we assume that $k<\ell$.

  \noindent {\bf Case 1.} Assume that $\j=\empwd$. Let $m=|\i|$. Then $m+k=\ell$ and $m\geq 1$.
  Notice that $$C^\ell=R_{\ell}(T_{i_0})\cup L_{\tau(\ell)}(T_{i_0+1})
  =\left(\bigcup_{j=n-\beta+1}^n T_{i_0[n]^\ell j}\right) \cup \left(\bigcup_{j=1}^\alpha T_{(i_0+1)[1]^{\tau(\ell)} j}\right).$$
  From $C_\i^k\subset T_\i$ and $\gt(\ell)\geq \ell>m$, we know that $C_\i^k\cap C^\ell=\emptyset$ if $\i\not\in \{i_0[n]^{m-1}, (i_0+1)[1]^{m-1}\}$.

  In case that $\i=i_0[n]^{m-1}$, we have
  $$C_\i^k\subset \gY_{i_0[n]^{m-1}}(C^1), \quad R_\ell(T_{i_0})=\gY_{i_0[n]^{m-1}}\left( \bigcup_{j=n-\beta+1}^n T_{[n]^{\ell-m+1}j}\right).$$
  Notice that
  \begin{align*}
    &\max C^1=\gY_{(i_0+1)[1]^{\tau(1)}\ga}(1)\leq \gY_{n1\ga}(1)<\gY_{n1}(1), \\
    &\min T_{[n]^{\ell-m+1}(n-\gb+1)}=\gY_{[n]^{\ell-m+1}(n-\gb+1)}(0)
    >\gY_{n(n-\gb+1)}(0)>\gY_{n(n-\gb)}(1)>\max C^1.
  \end{align*}
  We have $C_\i^k\cap C^\ell=C_\i^k\cap R_\ell(T_{i_0})=\emptyset$.

  In case that $\i=(i_0+1)[1]^{m-1}$, we have
  \begin{equation*}
    C_\i^k\subset \gY_{(i_0+1)[1]^{m-1}}(C^1), \quad L_{\tau(\ell)}(T_{i_0+1})=\gY_{(i_0+1)[1]^{m-1}}\left( \bigcup_{j=1}^\ga T_{[1]^{\tau(\ell)-m+1}j}\right).
  \end{equation*}
  Notice that
  \begin{align*}
    &\min C^1=\gY_{i_0 n(n-\gb+1)}(0)> \gY_{1n(n-\gb)}(1)>\gY_{1(n-1)}(1), \\
    &\max T_{[1]^{\tau(\ell)-m+1}\ga}<\gY_{1\ga}(1)\leq \gY_{1(n-1)}(1)<\min C^1,
  \end{align*}
  where we use $\tau(\ell)-m\geq\ell-m=k\geq 1$. We have $C_\i^k\cap C^\ell=C_\i^k\cap L_{\tau(\ell)}(T_{i_0+1})=\emptyset$.

  \noindent{\bf Case 2.} Assume that $\j\not=\empwd$.
  Let $\u\in\Sigma_n^*\cup\{\empwd\}$
  be the word with the maximal length which satisfies $\i=\u\i'$ and
  $\j=\u\j'$ for some $\i',\j'\in \Sigma_n^*\cup \{\empwd\}$.

  Suppose that $\u\not=\j$, then $\i',\j'$ are all in $\Sigma_n^*$
  with $\i'(1)\not=\j'(1)$, where $\i'(1)$ and $\j'(1)$ is the first letter of $\i'$ and $\j'$, respectively. Using $C_{\i'}^k\subset T_{\i'(1)}$ and $C_{\j'}^\ell\subset T_{\j'(1)}$, it is easy to see that $C_{\i'}^k\cap C_{\j'}^\ell=\emptyset$ so that $C_\i^k\cap C_\j^\ell=\gY_\u(C_{\i'}^k\cap C_{\j'}^\ell)=\emptyset.$

  Suppose that $\u=\j$. Using the result of Case 1, we have $C_\i^k\cap C_\j^\ell=\gY_\j(C_{\i'}^k\cap C^\ell)=\emptyset.$
\eproof

\medskip

\begin{lemma}\label{lem:C-ik-nonarchi}
  For any $A\in \sC_u$ and $B\in \sC_v$ with $u>v$. We have either $A\cap B=\emptyset$ or $A\subset B$.
\end{lemma}
\Proof Suppose that $A=C_\i^k$ and $B=C_\j^\ell$.

If $\i=\j=\empwd$, then the lemma holds in this case since $C^k\subset C^\ell$ for $k>\ell$. If $\i=\empwd$ and $\j\in\gS_n^*$, then from Lemma~\ref{lem:C-ik-disjoint}, we have $C^k\cap C_\j^\ell\subset C^{|\j|+\ell}\cap C_\j^\ell=\emptyset$ so that the lemma also holds in this case. Thus, we can assume that $\i\in\gS_n^*$ in the following.

Given $\i\in\gS_n^*$. It is easy to check that we must have either $C_\i^1\cap R_\ell(T_{i_0})=\emptyset$ or $C_\i^1\subset R_\ell(T_{i_0})$, while $C_\i^1\subset R_\ell(T_{i_0})$ if and only if one of the followings happens:\; (1).\; $\i=i_0[n]^\ell$ and $i_0\geq n-\gb+1$, or (2).\; $\i=i_0[n]^\ell j\u$ for some $j\in\{n-\gb+1,\ldots,n\}$ and $\u\in\gS_n^*\cup \{\empwd\}$. Similarly, we must have either $C_\i^1\cap L_{\tau(\ell)}(T_{i_0+1})=\emptyset$ or $C_\i^1\subset L_{\tau(\ell)}(T_{i_0+1})$, while $C_\i^1\subset L_{\tau(\ell)}(T_{i_0+1})$ if and only if one of the followings happens:\; (1).\; $\i=(i_0+1)[1]^{\tau(\ell)}$ and $i_0\leq \ga-1$, or (2).\; $\i=(i_0+1)[1]^{\tau(\ell)} j\u$ for some $j\in\{1,\ldots,\ga\}$ and $\u\in\gS_n^*\cup \{\empwd\}$. It follows that we must have either $C_\i^1 \cap C^{\ell}=\emptyset$ or $C_\i^1 \subset C^{\ell}$.

\noindent{\bf Case 1.} Assume that $\j=\empwd$. Since $C_\i^k\subset C_\i^1$ for any $k\in\bZp$, we know from above that the lemma holds in this case.

\noindent{\bf Case 2.} Assume that $\j\not=\empwd$.
  Let $\u\in\Sigma_n^*\cup\{\empwd\}$
  be the word with the maximal length which satisfies $\i=\u\i'$ and
  $\j=\u\j'$ for some $\i',\j'\in \Sigma_n^*\cup \{\empwd\}$.

  Suppose that both $\i'$ and $\j'$ are in $\gS_n^*$. Using the same method in Case~2 in the proof of Lemma~\ref{lem:C-ik-disjoint}, we have $C_\i^k\cap C_\j^\ell=\gY_\u(C_{\i'}^k\cap C_{\j'}^\ell)=\emptyset.$

  Suppose that one of $\i'$ and $\j'$ equals $\empwd$. Using the above discussions, we can easily see that one of the followings must holds: $C_\i^k\cap C_\j^\ell=\emptyset$ or $C_\i^k\subset C_\j^k.$
\eproof

\medskip


  Let $E$ be a given subset of $\bR$ and $\sP$ a family of finitely many closed subsets of $E$.
  If $\bigcup_{A\in \sP} A =E$ and the union is disjoint, we call $\sP$
  a \emph{partition} of $E$ and define $\|\sP\|=\max_{A\in \sP}\diam
  A$. Let $\sA_1$ and $\sA_2$ be
  two partitions of $E$.  If for any $A\in \sA_1$, there exist $j\in \bZp$ and
  $A'_1,\ldots,A'_j \in \sA_2$ such that $A=\bigcup_{i=1}^j A'_i$, then
  $\sA_2$ is called a \emph{refinement} of $\sA_1$. Clearly,
  $\sA_2$ is a refinement of $\sA_1$ if and only if for each $B\in
  \sA_2$, there exists $B'\in \sA_1$ such that $B\subset B'$.

  Let $\{\sP_k\}$ be a sequence of partitions of $T$. $\{\sP_k\}$ is
  called \emph{hierarchical} if
  $\sP_{k+1}$ is a refinement of $\sP_k$ for any $k$. $\{\sP_k\}$ is
  called \emph{convergent} if it is hierarchical and $\lim_{k\to
  \infty}\|\sP_k\|=0$.

  Denote by $\card A$ the cardinality of $A$ for any set $A$.
  Given a bounded subset $B$ of $\bR$. We define $CH(B)$ to be the convex hull of $B$. Equivalently, $CH(B)$ is the minimal closed interval containing $B$.

  Let $A$ be a given compact subset of $\bR$. Let $\{A_i\}_{i=1}^k$ be a family of compact subsets of $A$ with the following properties: $A_i$ is $A$-separate for all $i$, $CH(A_i)\cap A=A_i$  for all $i$, and $CH(A_i)$ does not intersect $CH(A_j)$ for all distinct $i$ and $j$. We define $\sS$ to be the family of compact subsets of $A$ with the minimal cardinality such that the following two conditions hold: (1). $\bigcup_{B\in \sS} B=A$ and $CH(B)\cap CH(B')=\emptyset$ for all distinct $B,B'\in \sS$; (2). $A_i\in \sS$ for all $i$. We call $\sS$ the \emph{simple decomposition} of $A$ by $\{A_i\}_{i=1}^k$. Clearly, there exists a unique simple decomposition for given $A$ and $\{A_i\}_{i=1}^k$. Furthermore, it is obvious that  we have the following property by definition:
  \begin{equation}\label{eq:simp-deco-prop}
    \mbox{\emph{ If $B\subset A$ and $CH(B)\cap \bigcup_{i=1}^k A_i=\emptyset$, then there exists a unique $E\in\sS$, such that $B\subset E$.}}
  \end{equation}
  We call this the \emph{containing property} of the simple decomposition. For convenience, $\sS=\{A\}$ is defined to be the simple decomposition of $A$ by $\emptyset$. It is clear that the containing property still holds in this case.

  Given $k\in \bZp$ and a compact subset $F$ of $T$, we define $\sC_k(F)=\{A: A\in \sC_k \mbox{ and } A\subset F\}$. Notice that $A$ is $F$-separate and $CH(A)\cap F=A$  for all $A\in\sC_k(F)$. We define $\sS_k(F)$ to be the simple decomposition of $F$ by $\sC_k(F)$.

Now we inductively construct $\{\sS_k\}$ as follows. Define $\sS_1=\sS_1(T)$ and $\sS_{k+1}=\bigcup_{F\in \sS_k} \sS_{k+1}(F)$ for $k\geq 1$. Clearly, $\sS_1$ is the simple decomposition of $T$ by $\{C^1\}$, i.e. $\sS_1=\{[0,a]\cap T, C^1, [b, 1]\cap T\}$, where $a=\gY_{i_0n(n-\gb)}(1)$ and $b=\gY_{(i_0+1)[1]^{\tau(1)}(\ga+1)}(0)$.

\begin{lemma}\label{lem:SK-hierar}
$\{\sS_k\}$ is a hierarchical partition sequence of $T$ such that $\sC_k\subset \sS_k$ for all positive integers $k$.
\end{lemma}
\Proof By the definition of the simple decomposition, we know that $\sS_{k+1}$ is a refinement of $\sS_k$ for all $k$, and $\sS_k$ is a partition of $T$ for all $k$. Thus, in order to prove the lemma, it suffices to show that $\sC_k\subset \sS_k$ for all $k$. We will prove this by induction. Clearly, $\sC_1\subset \sS_1$.

Assume that $\sC_k\subset \sS_k$ for all $k\leq m$ for some given $m\in\bZp$.

\begin{claim}

For any $A\in \sC_{m+1} $, there exists $B\in \sS_m$ such that $A\subset B$.

\end{claim}

\Proof
 Since $A\subset T$ and $\sS_m$ is a partition of $T$, there exists $B\in \sS_m$ such that $A\cap B\not=\emptyset$. Denote this $B$ by $B_m$. Notice that by definition, $\sS_{k+1}$ is a refinement of $\sS_k$ for all $k$. Hence there exists (unique) $\{B_k\}_{k=1}^{m-1}$ such that $B_{k+1}\subset B_k$ and $B_k\in \sS_{k}$ for all $k=1,2,\ldots,m-1$. It follows from $A\cap B_m\not=\emptyset$ that $A\cap B_k\not=\emptyset$ for all $k=1,2,\ldots,m$.

  Suppose that $A\cap C^1\not=\emptyset$. Then by Lemma~\ref{lem:C-ik-nonarchi}, we have $A\subset C^1$. Since all sets in the family $\sS_1$ are disjoint, we have $B_1=C^1$ so that $A\subset B_1$.
 Suppose that $A\cap C^1=\emptyset$. Notice that $A\subset T$ and $\sS_1$ is the simple decomposition of $T$ by $\{C^1\}$. Using the containing property of the simple decomposition, there exists a unique $E\in \sS_1\setminus \{C^1\}$ such that $A\subset E$. By the same reason as above, we have $B_1=E$ so that $A\subset B_1$. Hence, we always have $A\subset B_1$. Thus the claim holds in case that $m=1$.

 Assume that $m\geq 2$.
 Suppose that there exists $B'\in \sC_2$ such that $A\cap B'\not=\emptyset$. Similarly as above, we have $A\subset B'$. By the inductive assumption, we have $\sC_2\subset \sS_2$. Since $\sS_2$ is a partition of $T$, we have $B_2=B'$ so that $A\subset B_2$. Suppose that $A\cap \bigcup_{F\in \sC_2} F=\emptyset$. Then we have $A\cap \bigcup_{F\in \sC_2(B_1)} F=\emptyset$. Since $A\subset B_1$, we can obtain from the containing property of the simple decomposition that there exists $B'\in \sS_2(B_1)\setminus \sC_2(B_1)\subset \sS_2$ such that $A\subset B'$. Similarly as above, we have $B_2=B'$ so that $A\subset B_2$. Repeating this process, we can see that $A\subset B_k$ for $k=1,2,\ldots,m$. Thus the claim also holds in this case.
 \eproof

   From the above claim, we know that for each $A\in \sC_{m+1}$, there exists $B\in \sS_m$ such that $A\subset B$. Thus, we have $A\in \sS_{m+1}(B)\subset \sS_{m+1}$. It follows that  $\sC_k\subset \sS_k$ for $k=m+1$. By induction, $\sC_k\subset \sS_k$ for all $k\in \bZp$.
\eproof

\medskip

By Lemma~\ref{lem:SK-hierar}, we can show that the following corollary holds.

\begin{corol}\label{corol:conv-part-seq}
 There exists a
convergent partition sequence $\{\sT_k\}$ of $T$ such that
$\sC_k\subset \sT_k$ for all positive integers $k$.
\end{corol}
\Proof Let $E$ be a compact subset of $T$.
An open interval $(a,b)$ is said to be a \emph{gap} of $E$ if $a,b\in E$ and $(a,b)\cap E=\emptyset$. We call $b-a$ the \emph{length} of the gap $(a,b)$. Let $\gd$ be a positive real number. We define $\sG(E,\gd)=\{(a,b):\, (a,b) \textrm{ is a gap of $E$ such that $b-a\geq \gd$}\}$. Define $\sJ(E,\gd)$ to be the family of all connected components of $CH(E)\setminus \bigcup_{F\in \sG(E,\gd)} F$. Define $\sP(E,\gd)=\{A\cap E:\, A\in \sJ(E,\gd)\}$. Then $\sP(E,\gd)$ is a partition of $E$. Furthermore, for all $F\in \sP(E,\gd)$, $CH(F)$ does not contain any gap of $T$ whose length greater than $\gd$.

Now, we define $\gd_k=\max\{\diam A:\, A\in \sC_k\}$ for $k\in \bZp$. Then the sequence $\{\gd_k\}_{k=1}^\infty$ is decreasing and $\lim_{k\to \infty} \gd_k=0$. Define
\begin{equation*}
  \sT_k=\bigcup_{E\in \sS_k} \sP(E,\gd_k), \quad \forall k.
\end{equation*}
Clearly, $\sT_k$ is also a hierarchical partition sequence of $T$ with $\sC_k\subset \sT_k$ for all $k$. From $C_\i^1 \subset I_\i$, we can see that for any $A\in \sI_k$ with $k\in \bZ^+$, there exists $B\in \sC_{k+1}$ such that $B\subset A$. Thus
$$\|\sT_{k+1}\| < 2\cdot \max\{ \diam A:\, A\in \sI_k\} + \gd_{k+1}$$
for all $k$ so that $\lim_{k\to \infty} \| T_k\|=0$. It follows that the corollary holds.
\eproof


\subsection{\emph{Martingales and the proof of Theorem~\ref{thm:necessary-cond}}}

  Assume that $f: T\to D$ is bi-Lipschitz, i.e., $f$ is bijective and there exist two
  positive constants $c,c^\prime$ such that
  \begin{equation}\label{eq:biLipFct-def}
    c|x-y| \leq |f(x)-f(y)| \leq c^\prime |x-y|,\quad \forall x,y\in T.
  \end{equation}
  Let $s$ be the common
  Hausdorff dimension of $T$ and $D$, i.e. $\dim_H T=\dim_H D=s$.

  By Lemmas~\ref{lem:finite} and \ref{lem:T-gl-separated},
  there exists an integer $n_0$ such that for any $\i\in \Sigma_n^* \cup
  \{\empwd\}$ and $k\in \bZp$, there exist $\j, \j_1,\ldots,\j_p\in
  \Sigma_n^*$ such that $D_{\j\j_1}, D_{\j\j_2},\ldots,D_{\j\j_p}$
  are disjoint and
  \begin{equation*}
    f(C_{\i}^k) = \bigcup_{r=1}^p D_{\j\j_r} \subset D_{\j},
  \end{equation*}
  where each $|\j_r|=n_0$. Furthermore, by Remark~\ref{rem:optimal-decomp}, we can require
  $D_\j$ to be the smallest cylinder containing $f(C_{\i}^k)$.
  We denote this  $\j$ by $\j(\i,k)$ and define $\gga_{\i,k}=\sum_{r=1}^p
  \gr_{\j_r}^s$. Then
  \begin{equation}\label{eq:HM-f-C-i-k}
    \sH^s(f(C_\i^k)) = \sH^s (D_{\j(\i,k)}) \cdot \gga_{\i,k}, \quad \mbox{and}
  \end{equation}
  \begin{equation}\label{eq:D-j-ik-subset}
    D_{\j(\i,k)}\subset D_{\j(\i',k')} \quad \textrm{if} \quad  C_{\i}^k\subset
    C_{\i'}^{k'}.
  \end{equation}
  Define
  \begin{equation}\label{eq:def-sM}
    \sM=\left\{ \sum_{\j\in \sA} \gr_{\j}^s |\, \sA \subset
    \{1,\ldots,n\}^{n_0}\right\}.
  \end{equation}
  Then $\gga_{\i,k}\in\sM$ for all $\i$ and $k$.

  \medskip

  Let $\{\sT_k\}$ be a convergent partition sequence of $T$ as defined in
  Corollary~\ref{corol:conv-part-seq}.  We define
  \begin{equation*}
    g_k(A)=\frac{\sH^s(f(A))}{\sH^s(A)}, \quad A\in \sT_k.
  \end{equation*}
  Notice that for any $A\in \sT_k$, we can decompose $A$ by $A=\bigcup_{i=1}^j A'_i$
  where $A'_1,\ldots,A'_j\in
  \sT_{k+1}$. Then,
  \begin{eqnarray}\label{eq:martingale}
    \sum_{i=1}^j \frac{\sH^s(A'_i)}{\sH^s(A)}g_{k+1}(A'_i)
    = \sum_{i=1}^j \frac{\sH^s(A'_i)}{\sH^s(A)}
    \frac{\sH^s(f(A'_i))}{\sH^s(A'_i)}
    = \sum_{i=1}^j
    \frac{\sH^s(f(A'_i))}{\sH^s(A)}=\frac{\sH^s(f(A))}{\sH^s(A)}=g_k(A).
  \end{eqnarray}
  Let $\sF_k$ be the sigma field generated by $\sT_k$, and define
  \begin{equation*}
    g_k(x)=g_k(A)
  \end{equation*}
  for $x\in A$ where $A\in \sT_k$. Then $g_k$ is a
  $\sF_k$-measurable function. By (\ref{eq:martingale}), we know
  that $(g_k,\sF_k)$ is a martingale. Furthermore, by (\ref{eq:biLipFct-def}), we
  have $c^s \leq g_k(x) \leq (c^\prime)^s$
  for any $k$ and any $x\in T$. So the martingale convergence
  theorem implies that
  \begin{equation}\label{eq:martingale-result}
      g_k(x) \to g(x) \quad\textrm{as $k\to\infty$,} \quad \textrm{ for $\sH^s$-almost all } x \textrm{ in } T,
  \end{equation}
  where $g$ is $\sF$-measurable, with $\sF$
  the sigma field generated by $\bigcup_{k=1}^\infty \sF_k$.

  Define $\gm_{i}=\rho_i^s$ for any $i=1,\ldots,n$.
  For any $\i=i_1\cdots i_j\in\Sigma_n^*$, we define
  $\gm_\i =\prod_{k=1}^j \gm_{i_k}$. Denote $\gm_L=\sum_{j=1}^\ga \gm_j$ and $\gm_R=\sum_{j=n-\gb+1}^n \gm_j$. We have the following lemma.

  \begin{lemma}\label{lem:main-lemma}
    If $D\sim T$, then there exist a constant $M>0$ and infinitely
    many $(p_1,p_2,q_1,q_2)\in (\bZp)^4$ with $q_1\not=q_2$ and $|p_2-p_1|, |q_2-q_1|\leq
    M$ such that
    \begin{equation}\label{eq:main-eq}
          \frac{\gm_{i_0}\gm_n^{q_2}\gm_R + \gm_{i_0+1}\gm_1^{p_2}\gm_L}{\gm_{i_0}\gm_n^{q_1}\gm_R +
    \gm_{i_0+1}\gm_1^{p_1}\gm_L} = \gm_1^{k_1} \gm_2^{k_2}
    \cdots \gm_n^{k_n},
    \end{equation}
    where $\{k_i\}_{i=1}^n$ are nonnegative integers with $\max_i k_i\leq
    M$.
  \end{lemma}

  \Proof Let $\gS_n^\infty=\{i_1 i_2\cdots i_m\cdots |\,
  i_m\in\{1,\ldots,n\} \mbox{ for all } m\}$ . For each $x\in T$, there exists $i_1\cdots
  i_m\cdots\in\Sigma_n^\infty$ such that $\{x\}=\bigcap_{m\geq 1} \gY_{i_1\cdots
  i_m}(T)$. We call $i_1\cdots i_m\cdots$ to be the address of $x$. We
  remark that the address of $x$ may be not unique. However, if we
  define $\wdt{T}$ to be the set of all points in $T$ with unique
  address, then $\sH^s(\wdt{T})=\sH^s(T)$ by the definition of $T$. For
  each $x\in\wdt{T}$ with address $i_1\cdots i_m\cdots$, we define
  $\gs(x)$ to be the point with address $i_2\cdots
  i_m\cdots$. It is easy to check that $\gs(x)\in \wdt{T}$ for all $x\in\wdt{T}$.

  Let $\wdt{\sF}$ be the sigma field generated by $\{T_\i\cap\wdt{T}: \,
  \i\in\gS_n^*\}$. Define $\gn(A)=\sH^s(A\cap \wdt{T}) / \sH^s(\wdt{T}), \forall A\in
  \wdt{\sF}$. Then
  $\gs:\, (\wdt{T},\wdt{\sF},\gn)\to (\wdt{T},\wdt{\sF},\gn)$ is measure
  preserving. Fix $p\geq \card \sM+1$ in the proof of the lemma,
  where $\sM$ is defined by (\ref{eq:def-sM}). Given $q\in \bZp$, by the Poincar\'{e}
  recurrence theorem, for $\gn$-almost all $x\in C^{pq}\cap \wdt{T}$,
  i.e. for $\sH^s$-almost all $x\in C^{pq}\cap \wdt{T}$,
  there is an integer sequence $0<n_1(x,q)<n_2(x,q)<\cdots$ such
  that $\gs^{n_i(x,q)}(x)\in C^{pq}\cap\wdt{T}$ for all $i$. Thus, from
  (\ref{eq:martingale-result}), we can pick a point $x_q\in C^{pq}\cap
  \wdt{T}$ with $\gs^{n_i(x_q,q)}(x_q)\in C^{pq}\cap\wdt{T}$ for each $i$ and
  $g_k(x_{q})\to g(x_{q})$ as $k\to \infty$.

  Let $i_1\cdots i_m\cdots$ be the address of $x_q$. Define $\i_k=i_1i_2\cdots
  i_{n_k(x_q,q)}$. Then
  $$x_q=\gY_{\i_k}(\gs^{n_k(x_q,q)}(x_q)) \in  \gY_{\i_k}(C^{pq})=C_{\i_k}^{pq}.$$
  Recall that $
    \sH^s(f(C_{\i_k}^{t})) = \sH^s (D_{\j(\i_k,t)}) \cdot
    \gga_{\i_k,t}
  $ for any $t\in \bZp$.
  For any $t,t'\in \bZp\cap[1,p]$ with $t'<t$,
  using $x_q\in C_{\i_k}^{pq}\subset C_{\i_k}^{p(q-1)+t}\subset
  C_{\i_k}^{p(q-1)+t'}$, we have
  \begin{equation*}
    \frac{g_{|\i_k|+p(q-1)+t}(x_q)}{g_{|\i_k|+p(q-1)+t'}(x_q)}=
    \frac{\sH^s(D_{\j(\i_k,p(q-1)+t)})}{\sH^s(D_{\j(\i_k,p(q-1)+t')})}
    \cdot \frac{\gga_{\i_k,p(q-1)+t}}{\gga_{\i_k,p(q-1)+t'}}
    \cdot \frac{\sH^s(C_{\i_k}^{p(q-1)+t'})}{\sH^s(C_{\i_k}^{p(q-1)+t})}.
  \end{equation*}
  By (\ref{eq:D-j-ik-subset}), $D_{\j(\i_k,p(q-1)+t)}\subset
  D_{\j(\i_k,p(q-1)+t')}$ so that there exists $\u(q,k,t,t')\in
  \Sigma_n^*$ such that
  \begin{equation*}
    \sH^s(D_{\j(\i_k,p(q-1)+t)}) = \gm_{\u(q,k,t,t')}
    \sH^s(D_{\j(\i_k,p(q-1)+t')}).
  \end{equation*}
  Meanwhile, by definition,
  \begin{equation*}
    \frac{\sH^s(C_{\i_k}^{p(q-1)+t'})}{\sH^s(C_{\i_k}^{p(q-1)+t})} =
    \frac{\gm_{i_0}\gm_n^{p(q-1)+t'}\gm_R + \gm_{i_0+1}\gm_1^{\gt(p(q-1)+t')}\gm_L}{\gm_{i_0}\gm_n^{p(q-1)+t}\gm_R +
    \gm_{i_0+1}\gm_1^{\gt(p(q-1)+t)}\gm_L}.
  \end{equation*}
  Thus
  \begin{align}\label{eq:g-frac}
    \frac{g_{|\i_k|+p(q-1)+t}(x_q)}{g_{|\i_k|+p(q-1)+t'}(x_q)}=
    &\gm_{\u(q,k,t,t')}\cdot \frac{\gga_{\i_k,p(q-1)+t}}{\gga_{\i_k,p(q-1)+t'}}\notag\\
     &\cdot        \frac{\gm_{i_0}\gm_n^{p(q-1)+t'}\gm_R + \gm_{i_0+1}\gm_1^{\gt(p(q-1)+t')}\gm_L}{\gm_{i_0}\gm_n^{p(q-1)+t}\gm_R +
    \gm_{i_0+1}\gm_1^{\gt(p(q-1)+t)}\gm_L} .
  \end{align}

  \begin{claim}
    $\{\gm_{\u(q,k,t,t')}\}_{q\geq 1, k\geq 1, 1\leq t'<t\leq p}$ can take only finitely many values.
  \end{claim}
  \Proof Notice that $f$ is bi-Lipschitz. Thus for all $q\geq 1, k\geq 1, 1\leq t'<t\leq p$ we have
  \begin{equation*}
    \sH^s(f(C_{\i_k}^{p(q-1)+t})) \asymp \sH^s(C_{\i_k}^{p(q-1)+t}) \asymp \sH^s(C_{\i_k}^{p(q-1)}) \asymp
    \sH^s(C_{\i_k}^{p(q-1)+t'}) \asymp \sH^s(f(C_{\i_k}^{p(q-1)+t'})).
  \end{equation*}
  On the other hand, by (\ref{eq:HM-f-C-i-k}),
  \begin{equation*}
    \sH^s(D_{\j(\i_k,p(q-1)+t)}) \asymp
    \sH^s(f(C_{\i_k}^{p(q-1)+t})) \quad \textrm{and}\quad \sH^s(D_{\j(i_k,p(q-1)+t')}) \asymp
    \sH^s(f(C_{\i_k}^{p(q-1)+t'}))
  \end{equation*}
  for all $q\geq 1, k\geq 1, 1\leq t'<t\leq p$. Thus
  \begin{equation*}
    \sH^s(D_{\j(\i_k,p(q-1)+t)}) \asymp \sH^s(D_{\j(\i_k,p(q-1)+t')})
  \end{equation*}
  so that $\gm_{\u(q,k,t,t')} \asymp 1$. The claim follows
  immediately.
  \eproof

  By the claim, for fixed $q\geq 1$, the right hand of
  (\ref{eq:g-frac}) can take only finitely many distinct values for
  $1\leq t'<t\leq p$ and $k\in \bZp$ so that
  $\frac{g_{|\i_k|+p(q-1)+t}(x_q)}{g_{|\i_k|+p(q-1)+t'}(x_q)}$ can
  take only finitely many distinct values. Hence, we can take $k$ large
  enough such that
  $\frac{g_{|\i_k|+p(q-1)+t}(x_q)}{g_{|\i_k|+p(q-1)+t'}(x_q)}$ is
  so close to $1$ that it equals $1$. Since $p\geq \card \sM+1$, we can take $t_q>t'_q$
  such that $\gga_{\i_k,p(q-1)+t_q}=\gga_{\i_k,p(q-1)+t'_q}$. Thus,
  \begin{equation*}
    \frac{\gm_{i_0}\gm_n^{p(q-1)+t_q}\gm_R + \gm_{i_0+1}\gm_1^{\gt(p(q-1)+t_q)}\gm_L}{\gm_{i_0}\gm_n^{p(q-1)+t'_q}\gm_R +
    \gm_{i_0+1}\gm_1^{\gt(p(q-1)+t'_q)}\gm_L} = \gm_1^{k_1(q)} \gm_2^{k_2(q)}
    \cdots \gm_n^{k_n(q)},
  \end{equation*}
  where $\{k_i(q)\}_{1\leq i\leq n, q\geq 1}$ are bounded nonnegative integers. Also, we
  have $|(p(q-1)+t_q) - (p(q-1)+t'_q)| \leq p-1$. From
  \begin{equation*}
    \gr_1^{\gt(p(q-1)+t_q)} \asymp \gr_n^{p(q-1)+t_q} \asymp
    \gr_n^{p(q-1)+t'_q} \asymp \gr_1^{\gt(p(q-1)+t'_q)}
  \end{equation*}
  with respect to the indices $p,q$ we know that $\{\gt(p(q-1)+t_q) - \gt(p(q-1)+t'_q)\}_{q\geq 1}$ are
  bounded.

  Define $p_1=\gt(p(q-1)+t'_q),   q_1=p(q-1)+t'_q, p_2=\gt(p(q-1)+t_q),   q_2=p(q-1)+t_q.$
  From $1\leq t'_q<t_q\leq p$, we have $q_1,q_2\in \bZp\cap [p(q-1)+1,pq]$
  with $q_1\not=q_2$. Since $q$ can be arbitrary chosen in $\bZp$, we finally obtain infinitely many solution of
  (\ref{eq:main-eq}). The lemma is proved. \eproof

  \medskip

  Now, we can prove Theorem~\ref{thm:necessary-cond}.

  \noindent{\bf Proof of Theorem~\ref{thm:necessary-cond}.} \;
    By Lemma~\ref{lem:main-lemma}, there exist infinitely many
    solutions $(p_1,p_2,q_1,q_2)$ of (\ref{eq:main-eq}) with $q_1\not=q_2$ such that
    $(p_2-p_1)$, $(q_2-q_1)$, $\{k_i\}_{i=1}^n$ are constants. It
    follows that there are infinitely many $(p_1,q_1)\in (\bZp)^2$
    and constants $\gm_{i_0}\gm_n^{q_2-q_1}\gm_R$,
    $\gm_{i_0+1}\gm_1^{p_2-p_1}\gm_L$, $\gm_{i_0}\gm_R$, $\gm_{i_0+1}\gm_L$, $\gm_1^{k_1}\cdots
    \gm_n^{k_n}$ such that the following equation holds:
    \begin{equation*}
          \frac{(\gm_{i_0}\gm_n^{q_2-q_1}\gm_R)\cdot \gm_n^{q_1} +
          (\gm_{i_0+1}\gm_1^{p_2-p_1}\gm_L)\cdot \gm_1^{p_1}}{(\gm_{i_0}\gm_R)\cdot \gm_n^{q_1} +
        (\gm_{i_0+1}\gm_L)\cdot \gm_1^{p_1}} = \gm_1^{k_1} \gm_2^{k_2}
      \cdots \gm_n^{k_n}.
    \end{equation*}

    Assume that $\frac{\gm_{i_0}\gm_n^{q_2-q_1}\gm_R}{\gm_{i_0}\gm_R} \not
    = \frac{\gm_{i_0+1}\gm_1^{p_2-p_1}\gm_L}{\gm_{i_0+1}\gm_L} $. Then there
    is a constant $\gd$ such that $\frac{\gm_1^{p_1}}{\gm_n^{q_1}}=\gd$ for
    infinitely many $(p_1,q_1)\in (\bZp)^2$. Take $(p_1,q_1), (p'_1, q'_1)\in
    (\bZp)^2$ such that $(p_1,q_1)\not=(p'_1,q'_1)$ and
    $\frac{\gm_1^{p_1}}{\gm_n^{q_1}}=\gd=\frac{\gm_1^{p'_1}}{\gm_n^{q'_1}}$.
    It follows that $\gm_1^{p_1-p'_1} = \gm_n^{q_1-q'_1}$ so that $\gr_1^{p_1-p'_1} =
    \gr_n^{q_1-q'_1}$ which implies that $\log \gr_1/\log \gr_n \in
    \bQ$.

    Assume that $\frac{\gm_{i_0}\gm_n^{q_2-q_1}\gm_R}{\gm_{i_0}\gm_R}
    = \frac{\gm_{i_0+1}\gm_1^{p_2-p_1}\gm_L}{\gm_{i_0+1}\gm_L} $. Then
    $\gm_n^{q_2-q_1}=\gm_1^{p_2-p_1}$ with $q_2\not=q_1$ so that $\log \gr_1/\log \gr_n \in
    \bQ$.
  \eproof

\subsection{Algebraic dependence necessary condition for $n=4$}

In this subsection, we assume that $n=4$, $\gr_1=\gr_4$, $\gS_T=\{2\}$ and $f:\, T\to D$ is bi-Lipschitz. Denote $s:=\dim_H T=\dim_H D$. Denote $\gm_i:=\gr_i^s$ for each $i$.

Define $E_k=T_{2[4]^{k}}\cup T_{3[1]^{k}}$ for all $k\geq 0$.
Define $\sE_1=\{T_1,T_4,E_0\}$ and
\begin{equation*}
  \sE_{k+1} = \gY_1(\sE_k) \cup \gY_4(\sE_k)
  \cup \left( \gY_2(\sE_k)\setminus\{T_{2[4]^k}\}\right)
  \cup \left( \gY_3(\sE_k)\setminus\{T_{3[1]^k}\}\right) \cup
  \{E_{k}\}
\end{equation*}
for all $k\geq 1$. For example, the sets in the class $\sE_2$ are:
\begin{equation*}
  T_{11}, T_{14}, \gY_1(E_0), \;\; T_{41}, T_{44}, \gY_4(E_0), \;\;
  T_{21}, \gY_2(E_0), \;\; T_{34}, \gY_3(E_0), \;\; E_1.
\end{equation*}

\begin{remark}
  Let $i_0=2$ and $\{C^k\}_{k=1}^\infty$ defined as in subsection~2.2. Then $C^k=T_{2[4]^{k+1}}\cup T_{3[1]^{k+1}} =E_{k+1}$ for $k\geq 1$. Also, it is clear that the sequence $\{\sE_k\}_{k=1}^\infty$ is a convergent partition of $T$. However, in this subsection, we will not use these facts and the martingale convergent theorem.
\end{remark}

It is clear that the
following lemma holds.
\begin{lemma}\label{lem:sB-T-separate}
  Let $\gl_0=\min\{\frac{\diam(A)}{d(A,T\setminus A)}:\, A\in \sE_1\}$. Then each set in the family
  $\bigcup_{k=1}^\infty \sE_k$ is $(T,\gl_0)$-separate. 
\end{lemma}

For any $A\in \sE_k$, we define
\begin{equation*}
  \wdt{g}_k(A) = \frac{\sH^s(f(A))}{\sH^s(A)}.
\end{equation*}
We also abuse the notation $\wdt{g}_k(x)=\wdt{g}_k(A)$ for $x\in A$.
Assume that $\{A_1,\ldots, A_j\}$ be a partition of $A$ in
$\sE_{k+1}$, i.e. $A=\bigcup_{i=1}^j A_i$, $A_i\in \sE_{k+1}$ for each
$i$, and the union is disjoint. Then it is clear that
\begin{equation*}
  \wdt{g}_k(A) = \sum_{i=1}^j \frac{\sH^s(A_i)}{\sH^s(A)} \wdt{g}_{k+1}
  (A_i).
\end{equation*}

\begin{lemma}\label{lem:wdt-g-finite}
  The set $\{\frac{\wdt{g}_{k+1}(x)}{\wdt{g}_k(x)}:\, x\in T, k\geq
  1\}$ is finite.
\end{lemma}
\Proof Notice that $\frac{\sH^s(E_{k+1})}{\sH^s(E_{k})}=\gm_1$ for all
$k$. By induction, we can easily see that
\begin{equation*}
  \left\{ \frac{\sH^s(A)}{\sH^s(B)}:
  A\in \sE_{k+1}, B\in \sE_k  \;\; \textrm{with} \;\; A\subset B  \right\}
  =\left\{\gm_1,\gm_2,\gm_3,\gm_2+\gm_3,
  \frac{\gm_1\gm_2}{\gm_2+\gm_3},
  \frac{\gm_1\gm_3}{\gm_2+\gm_3}\right\}
\end{equation*}
for all $k\geq 1$. On the other hand, using Lemma~\ref{lem:finite}, Remark~~\ref{rem:optimal-decomp} and
 Lemma~\ref{lem:sB-T-separate}, and using the bi-Lipschitz property of
$f$, we can obtain that the set
\begin{equation*}
  \left\{ \frac{\sH^s(f(A))}{\sH^s(f(B))}:
  A\in \sE_{k+1}, B\in \sE_k  \;\; \textrm{with}
  \;\; A\subset B, \;\; k\geq 1  \right\}
\end{equation*}
is finite.

Given $x\in T$ and $k\geq 1$. We assume that $x\in A\subset
B$ with $A\in \sE_{k+1}$ and $B\in \sE_k$. Then
\begin{equation*}
  \frac{\wdt{g}_{k+1}(x)}{\wdt{g}_k(x)}
  =\frac{\sH^s (f(A))}{\sH^s(f(B))}
  \cdot \frac{\sH^s(B)}{\sH^s(A)}.
\end{equation*}
By above discussions, we know that the lemma holds. \eproof

Now we have the following property by using
Lemma~\ref{lem:wdt-g-finite}. We remark that the proof of this
property is same as the proof of Lemma~4 in \cite{XiRu08}, which was
restated in \cite{RRW} for completeness (see the proof of Lemma~2.4
therein). Thus we omit the proof.
\begin{lemma}\label{lem:measure-preserving}
  There is a set $A_0$ in the family $\bigcup_{k=1}^\infty \sE_k$
  and a constant $\gd>0$, such that $\wdt{g}_k(x)=\gd$ for all $x\in
  A_0$ and $k\geq k_0$, where $A_0\in \sE_{k_0}$. 
\end{lemma}

By the lemma, the restriction of $f$ on $A_{0}$ is
measure-preserving up to a constant. Thus, if we choose $\i_0\in
\gS_4^*$ such that $\gY_{\i_0}(T)\subset A_0$, then the restriction
of $f$ on $\gY_{\i_0}(T)$ is also measure-preserving up to a
constant. Hence, without loss of generality, we assume that
$A_0=\gY_{\i_0}(T)$ in the sequel.

By Lemmas~\ref{lem:finite} and ~\ref{lem:sB-T-separate}, there
exists an integer $n_1$ such that for any $k\in\bZ^+$, there exist
$\j,\j_1,\ldots,\j_p\in\gS_4^*$ such that $D_{\j\j_1},
D_{\j\j_2},\ldots,D_{\j\j_p}$ are disjoint and
\begin{equation*}
  f(\gY_{\i_0}(E_k))=\bigcup_{r=1}^p D_{\j\j_r},
\end{equation*}
where each $|\j_r|=n_1$. Furthermore, by
Remark~\ref{rem:optimal-decomp}, we can require $D_\j$ to be the
smallest cylinder containing $f(\gY_{\i_0}(E_k))$. We denote this
$\j$ by $\j^\prime(k)$ and define $\gga^\prime_{k}=\sum_{r=1}^p
\gm_{\j_r}$. Define
  \begin{equation*}\label{eq:def-sMp}
    \sM^\prime=\left\{ \sum_{\j\in \sA} \gm_{\j} |\, \sA \subset
    \{1,\ldots,4\}^{n_1}\right\}.
  \end{equation*}

\medskip

\noindent{\bf Proof of Theorem~\ref{thm:branch-4}.} \;
Notice that $2\gm_1+\gm_2+\gm_3=1$. In order to prove the lemma, it
suffices to show that there exists a polynomial $P(x_1,x_2,x_3)$ with rational coefficients
such that $P(\gm_1,\gm_2,\gm_3)=0$ and $P((1-x_2-x_3)/2, x_2,x_3)$ is not
identically equal to $0$.

Define $\wdt{f}: T\to D$ by
\begin{equation*}
  \wdt{f}(x) = f(\gY_{\i_0}(x)).
\end{equation*}
Let $x^*$ be the unique point in the set $T_2\cap T_3$. Assume that
$\{\wdt{f}(x^*)\}=D_{t_1t_2\cdots}$.

Given $\i=i_1i_2\cdots i_k\in \gS_4^*$ and $\ell\in\{1,2,3,4\}$, we
define $N_\ell(\i)$ to be the cardinality of $\{j:\, i_j=\ell, 1\leq
j\leq k\}$.

\noindent \textbf{Case 1.}~~Assume that there are infinitely many
$k$ such that $t_k=2$. Then $\lim_{q\to \infty}
N_2(\j^\prime(q))=\infty$. Notice that $\card M'<+\infty$.
Thus, we can choose $q_1,q_2$ with $1\leq
q_1<q_2$ such that
$\gga^\prime_{q_1}=\gga^\prime_{q_2}$ and
\begin{equation}\label{ineq:N2-jprime}
  N_2(\j^\prime(q_2))>N_2(\j^\prime(q_1)).
\end{equation}
From $\gga^\prime_{q_1}=\gga^\prime_{q_2}$ and
Lemma~\ref{lem:measure-preserving}, we have
\begin{equation*}
  \frac{\sH^s(\gY_{\i_0}(E_{q_2}))}{\sH^s(\gY_{\i_0}(E_{q_1}))}
  =\frac{\sH^s(D_{\j^\prime(q_2)})}{\sH^s(D_{\j^\prime (q_1)})}.
\end{equation*}
Thus, by (\ref{ineq:N2-jprime}),
\begin{equation*}
  \gm_1^{q_2-q_1}=\gm_1^{k_1}\gm_2^{k_2}\gm_3^{k_3}
\end{equation*}
with $k_2\in\bZ^+$ and $k_1,k_3\in \bN$. Since the above equality
does not hold if we plug in $\gm_1=1/2,\gm_2=\gm_3=0$, we know that
$\gm_2$ and $\gm_3$ are algebraic dependent.

\noindent \textbf{Case 2.}~~Assume that there are infinitely many
$k$ such that $t_k=3$. Then using the same method as in Case 1, we
can obtain that $\gm_2$ and $\gm_3$ are algebraic dependent.

\noindent \textbf{Case 3.}~~Assume that there are only finitely many
$k$ such that $t_k\in\{2,3\}$. Then there exists $q_0\in \bZ^+$,
such that $N_2(\j^\prime(q))=N_2(\j^\prime(q_0))$ and
$N_3(\j^\prime(q))=N_3(\j^\prime(q_0))$ for all $q\geq q_0$.

By definition,
\begin{equation*}
  \frac{\sH^s(\wdt{f}(E_{q_0}))}{\sH^s(E_{q_0})}
  =\frac{\sH^s(D_{\j^\prime(q_0)})\gga^\prime_{q_0}}{
  (\gm_2+\gm_3)\gm_1^{q_0}\cdot\sH^s(T)}.
\end{equation*}
Substituting $\gm_3$ by $1-2\gm_1-\gm_2$, we know that $\gga^\prime_{q_0}$ is a polynomial of $\gm_1$ and $\gm_2$ with integral coefficients. By using Euclidean algorithm, it is easy to see that there exist polynomials $Q(\gm_1,\gm_2)$ and $R(\gm_2)$ with rational coefficients, such that $
 \gga^\prime_{q_0}  = (1-2\gm_1) Q(\gm_1,\gm_2)+R(\gm_2)$. Thus $\gga^\prime_{q_0}=(\gm_2+\gm_3)Q(\gm_1,\gm_2)+R(\gm_2)
$ so that
\begin{equation}\label{eq:Hs-ratio-Ap0}
  \frac{\sH^s(\wdt{f}(E_{q_0}))}{\sH^s(E_{q_0})}
  =\frac{\sH^s(D_{\j^\prime(q_0)})\big((\gm_2+\gm_3) Q(\gm_1,\gm_2)+R(\gm_2)\big)}{(\gm_2+\gm_3)\gm_1^{q_0}
  \cdot\sH^s(T)}.
\end{equation}

From $E_{q_0}=T_{2[4]^{q_0}}\cup T_{3[1]^{q_0}}$, we know that
$\gY_{\i_0}(T_{2[4]^{q_0}1})\subset \gY_{\i_0}(E_{q_0})$. Thus the
smallest cylinder of $D$ containing $f(\gY_{\i_0}(T_{2[4]^{q_0}1}))$
is a subset of the smallest cylinder of $D$ containing
$f(\gY_{\i_0}(E_{q_0}))$. It follows that there exists a polynomial
$P$ of $\gm_1,\gm_2,\gm_3$ with integral coefficients, such that
\begin{equation}\label{eq:Hs-ratio-2-4p0-1}
  \frac{\sH^s(\wdt{f}(T_{2[4]^{q_0}1}))}{\sH^s(T_{2[4]^{q_0}1})}
  =\frac{\sH^s(D_{\j^\prime(q_0)})\cdot
  P(\gm_1,\gm_2,\gm_3)}{\gm_1^{q_0+1}\gm_2\cdot\sH^s(T)}.
\end{equation}

By Lemma~\ref{lem:measure-preserving} and the definition of
$\wdt{f}$, the right hand sides of (\ref{eq:Hs-ratio-Ap0}) and
(\ref{eq:Hs-ratio-2-4p0-1}) are equal so that
\begin{equation}\label{eq:gga-q0-m1m2m3}
  \big((\gm_2+\gm_3) Q(\gm_1,\gm_2)+R(\gm_2)\big)\gm_1\gm_2 = (\gm_2+\gm_3) P(\gm_1,\gm_2,\gm_3).
\end{equation}

\noindent \textbf{Case 3.1.}~~Assume that the polynomial $R(x)$ is not identically equal to $0$. Then there exists $a\not=0$ such that $R(a)\not=0$. Notice that (\ref{eq:gga-q0-m1m2m3}) does not hold if we plug in $\gm_1=1/2$, $\gm_2=a$, $\gm_3=-a$. Thus $\gm_2$ and $\gm_3$ are algebraic dependent.

\noindent \textbf{Case 3.2.}~~Assume that the polynomial $R(x)$ is identically equal to $0$.
Then from (\ref{eq:Hs-ratio-Ap0}), we have
\begin{equation}\label{eq:Hs-ratio-polynomial}
  \frac{\sH^s(\wdt{f}(E_{q_0}))}{\sH^s(E_{q_0})}
  =\frac{\sH^s(D_{\j^\prime}(q_0))\cdot
  Q(\gm_1,\gm_2)}{\gm_1^{q_0}\cdot \sH^s(T)}.
\end{equation}

Since there are only finitely many $k$ such that $t_k=\{2,3\}$, we can
choose $q_3>q_0$ such that in the right hand side of
\begin{equation*}
  f(\gY_{\i_0}(E_{q_3}))=\bigcup_{r=1}^p D_{\j\j_r},
\end{equation*}
where $\j$, $\j_1,\ldots,\j_p$ have same meaning as above, the
cylinder $D_{\j\j_r}$ containing $\gY_{\i_0}(x^*)$ satisfies
$N_2(\j_r)=N_3(\j_r)=0$. It follows that $\gga^\prime_{q_3}\not=0$
if we plug in $\gm_1=1/2, \gm_2=\gm_3=0$.

Notice that
\begin{equation*}
  \frac{\sH^s(\wdt{f}(E_{q_3}))}{\sH^s(E_{q_3})}
  =\frac{\sH^s(D_{\j^\prime (q_3)})\cdot
  \gga^\prime_{q_3}}{(\gm_2+\gm_3)\gm_1^{q_3}\cdot \sH^s(T)}.
\end{equation*}
By (\ref{eq:Hs-ratio-polynomial}) and using
Lemma~\ref{lem:measure-preserving}, we have
\begin{equation*}
  (\gm_2+\gm_3)\gm_1^{q_3-q_0}Q(\gm_1,\gm_2)
  = \frac{\sH^s(D_{\j^\prime(q_3)})}{\sH^s(D_{\j^\prime(q_0)})}\cdot
  \gga^\prime_{q_3}.
\end{equation*}
It follows from $N_2(\j^\prime(q_3))=N_2(\j^\prime(q_0))$ and
$N_3(\j^\prime(q_3))=N_3(\j^\prime(q_0))$ that
$\frac{\sH^s(D_{\j^\prime(q_3)})}{\sH^s(D_{\j^\prime(q_0)})}
=\gm_1^{k}$
for some $k\in \bN$. Thus the above equality does not hold if we plug in $\gm_1=1/2, \gm_2=\gm_3=0$ so that $\gm_2$ and $\gm_3$ are
algebraic dependent. \eproof

\begin{example}
  Let $n=4$ and $\gS_T=\{2\}$. Let $\gm_2=e/4$, $\gm_3=1/4$ and $\gm_1=\gm_4=(1-\gm_2-\gm_3)/2$, where $e$ is the Euler constant. Let $\gr_i=\gm_i^2$ for all $i$. Then $\dim_H D=\dim_H T=1/2$. By Theorem~\ref{thm:branch-4}, $D\not\sim T$.
\end{example}


\section{Sufficient condition for $D\sim T$}

\subsection{\emph{Graph-directed sets corresponding to $D$ and $T$}}


%


Now we recall the notion of \emph{graph-directed set}. Let
$G=(V,\gGa)$ be a directed graph and $d$ a positive integer. Suppose for each edge $e\in \gGa$,
there is a corresponding similarity $S_e:\, \bR^d\to \bR^d$ with
ratio $r_e$. Assume that for each vertex $i\in V$, there exists an
edge starting from $i$, and assume that $r_{e_1}\cdots r_{e_k}<1$
for any cycle $e_1\cdots e_k$. Then there exists a unique family
$\{K_i\}_{i\in V}$ of compact subsets of $\bR^d$ such that for any
$i\in V$,
\begin{equation}\label{eq:GDS-def}
  K_i=\bigcup_{j\in V} \bigcup_{e\in \sE_{ij}} S_e(K_j),
\end{equation}
where $\sE_{ij}$ is the set of edge starting from $i$ and ending at
$j$. In particular, if the union in (\ref{eq:GDS-def}) is disjoint
for any $i$, we call $\{K_i\}_{i\in V}$ are dust-like graph-directed
sets on $(V,\gGa)$. For details on graph-directed sets, please see
\cite{MaWi88,RW98}.

Similarly as Theorem~2.1 in \cite{RRX06}, we have the following
lemma which was also pointed out in \cite{XiRu07}.

\begin{lemma}\label{lem:gds-strong}
  Suppose that $\{K_i\}_{i\in V}$ and $\{K_i^\prime\}_{i\in V}$ are dust-like
  graph-directed sets on $(V,\gGa)$ satisfying (\ref{eq:GDS-def})
  and $K_i^\prime=\bigcup_{j\in V} \bigcup_{e\in \sE_{ij}}
  S_e^\prime(K_j^\prime)$. If
  similarities $S_e$ and $S_e^\prime$ have the same ratio for each
  $e\in\gGa$ , then $K_i\sim
  K_i^\prime$ for each $i\in V$.
\end{lemma}

\begin{defi}
Assume that $\mathcal{K}=\{K_i\}_{i=1}^m$ and
$\mathcal{K}^\prime=\{K_i^\prime\}_{i=1}^m$ are two families of
compact subsets of $\bR^d$, where $m\geq 2$ is a given positive
integer. We say that two compact subsets $A$ and $B$ of $\bR^d$ have
same dust-like decomposition w.r.t. $\mathcal{K}$ and
$\mathcal{K}^\prime$, if there exist a positive integer $t\geq 2$
and positive integers $j_1,j_2,\ldots,j_t\in\{1,\ldots,m\}$, such
that
\begin{equation*}
   A=\bigcup_{i=1}^t S_{i}(K_{j_i}), \quad \mbox{and}\quad
   B=\bigcup_{i=1}^t S_{i}^\prime (K_{j_i}^\prime),
\end{equation*}
where the above two unions are disjoint, while for each $i$, $S_{i}$
and $S_{i}^\prime$ are similarities from $\bR^d$ to $\bR^d$ with
same ratio.
\end{defi}

Clearly, the following lemma is a weak version of
Lemma~\ref{lem:gds-strong}.

\begin{lemma}\label{lem:gds-weak}
  Suppose that $\mathcal{K}=\{K_i\}_{i=1}^m$ and
  $\mathcal{K}^\prime=\{K_i^\prime\}_{i=1}^m$ are two families of
  compact subsets of $\bR^d$. If for each $1\leq i\leq m$, $K_i$ and
  $K_i^\prime$ have same dust-like decomposition w.r.t.
  $\mathcal{K}$ and $\mathcal{K}^\prime$, then $K_i\sim
  K_i^\prime$ for each $i=1,\ldots,m$. 
\end{lemma}

Given $\i\in \gS_n^*$ and a nonnegative integer $k$.
We recall that
$L_k(T_\i)=\bigcup_{j=1}^{\ga} T_{\i[1]^k j}$ and
  $R_k(T_\i)=\bigcup_{j=n-\gb+1}^n T_{\i [n]^k j}$ as defined in (\ref{eq:Lk-Rk-def}).
Now we define
\begin{equation*}
  L_k(D_\i)=\bigcup_{j=1}^\ga D_{\i [1]^k j},
  \quad R_k(D_\i)=\bigcup_{j=n-\gb+1}^n
  D_{\i [n]^k  j}.
\end{equation*}

In the rest of this section, we will always assume that $\log
\gr_1/\log \gr_n \in \bQ$ and every touching letter is substitutable.
Thus, for any $i\in \gS_T$, there exist $\j_i\in\gS_n^*$ and
$k_i,k_i^\prime \in\bN$, such that one of the following holds:

(1)\; $\diam L_{k_i}(T_{i+1})=\diam L_{k_i^\prime}(T_{i\j_i})$  and the last letter of
$\j_i$ does not belong to $\{1\}\cup(\gS_T+1)$,

(2)\; $\diam R_{k_i}(T_{i})=\diam R_{k_i^\prime}(T_{(i+1)\j_i})$ and the last letter of
$\j_i$ does not belong to $\{n\}\cup\gS_T$.

From $\log \gr_1/\log \gr_n\in \bQ$, there exist $p,q\in \bZp$ such that $\gr_1^p=\gr_n^q$. Notice that we can choose $p$ and $q$ large enough such that $p,q>\max\{ k_i^\prime +|\j_i|: \; i\in \gS_T\}$. As a result, we have
\begin{equation}\label{eq:pq-left-restriction}
  L_{k_i^\prime}(T_{i[n]^{2q}\j_i}) \cap R_{3q}(T_i)
  =L_{k_i^\prime}(D_{i[n]^{2q}\j_i}) \cap R_{3q}(D_i) =\emptyset \quad
  \mbox{and} \quad
  L_{k_i^\prime}(D_{i\j_i}) \cap R_{q}(D_i)=\emptyset
\end{equation}
for any left substitutable touching letter $i$, and
\begin{equation}\label{eq:pq-right-restriction}
  L_{3p}(T_{i+1}) \cap
  R_{k_i^\prime}(T_{(i+1)[1]^{2p}\j_i})
  =L_{3p}(D_{i+1}) \cap
  R_{k_i^\prime}(D_{(i+1)[1]^{2p}\j_i})
  =\emptyset \quad
  \mbox{and}\quad
  L_{p}(D_{i+1}) \cap
  R_{k_i^\prime}(D_{(i+1)\j_i})=\emptyset
\end{equation}
for any right substitutable touching letter $i$.
We remark that this restriction is useful in the definition of
$\Dfour_i$ and in the proof of
Lemma~\ref{lem:gds-TD-34}. We will fix $p,q$ in this section.

Given a positive integer $m$. Let $\sJ_m=\bigcup_{|\i|=m} I_\i$ and
$J_{m,1},\ldots,J_{m,c_m}$ are the connected components of $\sJ_m$,
spaced from left to right. We define $\gL_i=\{j|\, I_j\subset
J_{1,i}\}$ for $i=1,\ldots,c_1$. Define
\begin{equation*}
  \Tone_i = \bigcup_{j\in \gL_i} T_j, \qquad \Done_i=\bigcup_{j\in
  \gL_i} D_j, \quad \forall i=1,\ldots,c_1.
\end{equation*}
Then for any $\i\in\gS_n^*$ and $k\in \bN$, we have
\begin{eqnarray*}
  L_k(T_\i)= \gY_{\i [1]^k} (\Tone_1), \qquad R_k(T_\i)= \gY_{\i [n]^k}
  (\Tone_{c_1}), \\
  L_k(D_\i)= \gF_{\i [1]^k} (\Done_1), \qquad R_k(D_\i)= \gF_{\i [n]^k}
  (\Done_{c_1}).
\end{eqnarray*}
For any touching letter $i$, we define
\begin{align*}
  &\Ttwo_i=R_0(T_i)\cup L_0(T_{i+1}), \quad &&\Dtwo_i=R_0(D_i)\cup
  L_0(D_{i+1}),\\
  &\Tthree_i=R_q(T_i)\cup L_p(T_{i+1}), &&\Dthree_i=R_q(D_i)\cup
  L_p(D_{i+1}).
\end{align*}
We remark that $R_0(T_i)=\gY_i(\Tone_{c_1})$ and
$L_0(T_{i+1})=\gY_{i+1}(\Tone_1)$. Furthermore, for any touching
letter $i$, if $i$ is left substitutable, we define
\begin{equation*}
    \Tfour_i=R_q(T_i)\cup L_{k_i}(T_{i+1}),
    \qquad  \Dfour_i=L_{k_i^\prime}(D_{i\j_i})\cup
    R_q(D_i).
\end{equation*}
Otherwise, we define
\begin{equation*}
    \Tfour_i=R_{k_i}(T_i)\cup L_p(T_{i+1}),
    \qquad \Dfour_i=L_p(D_{i+1})\cup
    R_{k_i^\prime}(D_{(i+1)\j_i}).
\end{equation*}

Define $\sT=\{T\}\cup \{\Tone_i:\,
i=1,\cdots,c_1\}\cup\{T_i^{(j)}:\, i\in\gS_T, j=2,3,4\}$ and
$\sD=\{D\}\cup \{\Done_i:\, i=1,\cdots,c_1\}\cup\{D_i^{(j)}:\,
i\in\gS_T, j=2,3,4\}$. We will show that each corresponding pair in
$\sT$ and $\sD$ have same dust-like decomposition w.r.t. $\sT$ and
$\sD$.

\subsection{\emph{The family $\sT^*$}}

Given $\i\in \gS_n^*$. We say that $\i\in\gS^*_L$ if there exist
$\i^\prime\in\gS_n^*\cup\{\empwd\}$, $j\in\gS_T+1$ and $k\in\bN$
such that $\i=\i^\prime j [1]^k$. Similarly, we say that
$\i\in\gS^*_R$ if there exist $\i^\prime\in\gS_n^*\cup\{\empwd\}$,
$j\in\gS_T$ and $k\in\bN$ such that $\i=\i^\prime j [n]^k$.  The
following lemma is easy to check.
\begin{lemma}\label{lem:gY-i-Tonej-sep}
  Given $\i\in\gS_n^*$ and $j\in\{1,\ldots,c_1\}$. Then
  $\gY_\i(\Tone_j)$ is $T$-separate if and only if one of the following
  conditions holds: \; (1).\; $2\leq j\leq c_1-1$;\;\;
  (2).\; $j=1$ and $\i\not\in\gS_L^*$;\;\;
  (3).\; $j=c_1$ and $\i\not\in\gS_R^*$. 
\end{lemma}


Define
\begin{equation*}
  \sT^{(1)}=\{\gY_\i(\Tone_j):\, \i\in\gS_n^*\cup\{\empwd\}
  \textrm{ and }
  j\in\{1,\ldots,c_1\} \textrm{ such that $\gY_\i(\Tone_j)$ is
  $T$-separate}\},
\end{equation*}
\begin{equation*}
  \sT^{(k)}=\{\gY_\i(T_j^{(k)}):\, \i\in\gS_n^*\cup\{\empwd\}
  \textrm{ and }  j\in\gS_T\}, \quad k=2,3.
\end{equation*}
It is clear that all sets in $\sT^{(2)}$ and $\sT^{(3)}$ are $T$-separate.
Define
$$\sT^*=\{A|\, A \mbox{ is a disjoint union of finitely many
($\geq 2$) sets in the
class } \sT^{(1)} \cup \sT^{(2)}\cup \sT^{(3)}\}.$$

\begin{remark}\label{rmk:sT-map}
  Assume that $A\in \sTOne$ with $A\subset (0,1)$. Then it is easy to check that $\gY_\i(A)\in \sTOne$ for any $\i\in \gS_n^*$. It follows that $\gY_\i(B)\in\sT^*$ for all $B\in\sT^*$ with $B\subset (0,1)$ and all $\i\in\gS_n^*$.
\end{remark}

Let $\gS_n^\infty=\{i_1 i_2\cdots i_m\cdots |\,
i_m\in\{1,\ldots,n\} \mbox{ for all } m\}$ as defined in the proof of Lemma~\ref{lem:main-lemma}. Given $\i=i_1\cdots
i_m\cdots \in \gS_n^\infty$, there exists  a unique point $x\in T$
such that
\begin{equation*}
  \{x\} = \bigcap_{m=1}^\infty \gY_{i_1\cdots i_m} ([0,1]).
\end{equation*}
We denote this unique $x$ by $\gp_T(\i)$. Then $\gp_T:\,
\gS_n^\infty \to T$ is a surjection. Similarly, we can define
$\gp_D:\, \gS_n^\infty \to D$ by
\begin{equation*}
  \{\gp_D(\i)\} = \bigcap_{m=1}^\infty \gF_{i_1\cdots i_m} ([0,1]),
  \quad \forall\, \i=i_1\cdots i_m \cdots \in \gS_n^\infty.
\end{equation*}
Since $D$ is dust-like, $\gp_D$ is a bijection.

By definition of $\gp_T$ and $\gp_D$, it is easy to check that
\begin{equation}\label{eq:gpD-gpT-relation}
  \gp_D\circ \gp_T^{-1} (\gY_\i(T_j^{(k)}))=\gF_\i(D_j^{(k)}),
  \quad \forall  \gY_\i(T_j^{(k)})\in \sT^{(k)}, \; k=1,2,3.
\end{equation}
Using this fact, we have the following lemma.

\begin{lemma}\label{lem:gpDgpTInv-decom}
  Let $A\in\sT^*$. Then $A$ and $\gp_D\circ \gp_T^{-1}(A)$ have
  same dust-like decomposition w.r.t. $\sT$ and $\sD$.
\end{lemma}
\Proof
  Assume that $A=\bigcup_{i=1}^m A_i$, where $m\geq 2$, $A_i\in \sT^{(1)} \cup \sT^{(2)}\cup \sT^{(3)}$ and the union
  is disjoint. Then
  \begin{equation}\label{eq:gdgp-decomp-sT}
    \gp_D\circ \gp_T^{-1}(A) = \bigcup_{i=1}^m \gp_D\circ \gp_T^{-1}
    (A_i).
  \end{equation}
  From the union in $\bigcup_{i=1}^m A_i$ is disjoint, we can see that the union in $\bigcup_{i=1}^m \gp_T^{-1}(A_i)$ is disjoint. Since $\gp_D$ is a bijection, we know that the union in (\ref{eq:gdgp-decomp-sT}) is also disjoint. From
  (\ref{eq:gpD-gpT-relation}), we can see that the lemma holds.
\eproof

\subsection{\emph{Graph-directed decomposition of $\sT$ and $\sD$ and the proof of sufficient condition}}

The following lemma is easy to show.
\begin{lemma}\label{lem:gds-TD-01}
  The following pairs have same dust-like decomposition w.r.t. $\sT$
  and $\sD$.
  \begin{enumerate}
    \item $T$ and $D$;
    \item $\Tone_i$ and $\Done_i$ for $i=1,\cdots,c_1$;
  \end{enumerate}
\end{lemma}
\Proof (i)~~ Clearly, $T$ and $D$ can be decomposed to following
disjoint unions.
\begin{equation*}
  T=\bigcup_{i=1}^{c_1} \Tone_i, \qquad D=\bigcup_{i=1}^{c_1} \Done_i.
\end{equation*}

(ii)~~ Given $i=1,\ldots,c_1$. Notice that
\begin{equation}\label{eq:Tone-decomp}
  \Tone_i=\bigcup_{j\in \gL_i} T_j =\bigcup_{j\in \gL_i}\gY_j(T)
  = \bigcup_{j\in\gL_i}
  \bigcup_{k=1}^{c_1} \gY_j(\Tone_k).
\end{equation}

Let $b(i)$ and $e(i)$ be the minimal and maximal element in $\gL_i$,
respectively. If $b(i)=e(i)$, then
\begin{equation*}
  \Tone_i=\bigcup_{k=1}^{c_1} \gY_{b(i)} (\Tone_k) \in \sT^*,
\end{equation*}
since each $\gY_{b(i)}(\Tone_k)$ is $T$-separate in this case. If
$b(i)<e(i)$, then for any $b(i)\leq j<e(i)$,
\begin{equation*}
  \gY_j(\Tone_{c_1})\cup \gY_{j+1}(\Tone_1) = \Ttwo_j,
\end{equation*}
and other $\gY_j(\Tone_k)$ in (\ref{eq:Tone-decomp}) are
$T$-separate so that they belong to $\sTOne$. Thus $\Tone_i$ also
belongs to $\sT^*$ in this case. By
Lemma~\ref{lem:gpDgpTInv-decom}, $\Tone_i$ and
$\Done_i=\gp_D\circ\gp_T^{-1}(\Tone_i)$ have same dust-like
decomposition w.r.t. $\sT$ and $\sD$.
\eproof

\medskip

\begin{remark}\label{rmk:Tone-T}
It follows from the above lemma that $\sTOne\subset \sT^*$.
\end{remark}

\begin{lemma}\label{lem:LuLv-SC}
  Given $\i\in\gS_n^*$ and two nonnegative integers $u,v$ with $u<v$. The pairs $L_u(T_\i)\setminus L_v(T_\i)$ and $L_u(D_\i)\setminus L_v(D_\i)$,
  $R_u(T_\i)\setminus R_v(T_\i)$ and $R_u(D_\i)\setminus R_v(D_\i)$ have same dust-like decomposition w.r.t. $\sT$ and $\sD$, respectively.
\end{lemma}
\Proof Without loss of generality, we only show that
  $L_u(T_\i)\setminus L_v(T_\i)$ and $L_u(D_\i)\setminus L_v(D_\i)$ have same dust-like decomposition w.r.t. $\sT$ and $\sD$. Clearly,
  $$\gp_D\circ \gp_T^{-1}(L_u(T_\i)\setminus L_v(T_\i))
  =L_u(D_\i)\setminus L_v(D_\i).$$
  Thus, from Lemma~\ref{lem:gpDgpTInv-decom} and noticing that $L_u(T_\i)\setminus L_v(T_\i)=\bigcup_{k=u}^{v-1}
  \Big(L_k(T_\i)\setminus L_{k+1}(T_\i)\Big)$, it suffices to show
  that $L_k(T_\i)\setminus L_{k+1}(T_\i)\in \sT^*$ for all
  $k\in \bN$.

  Given $k\in\bN$. Assume that $1\not\in \gS_T$, i.e. $\ga=1$. Then
  \begin{equation*}
    L_k(T_\i)\setminus L_{k+1}(T_\i) = \bigcup_{j=2}^n T_{\i [1]^{k+1}j}
    =\bigcup_{j=2}^{c_1}\gY_{\i [1]^{k+1}}(\Tone_j).
  \end{equation*}
  Notice that $\gY_{\i [1]^{k+1}}(\Tone_{c_1})$ is $T$-separate in this
  case. By Remark~\ref{rmk:Tone-T}, $L_k(T_\i)\setminus L_{k+1}(T_\i)\in \sT^*$.

  Assume that $1\in\gS_T$, i.e. $\ga\geq 2$. Then
  $L_k(T_\i)\setminus L_{k+1}(T_\i)=\gY_{\i [1]^k}(A)$, where
  \begin{equation}\label{eq:A-decomp}
    A = \left(\bigcup_{j=\ga+1}^n T_{1j}\right) \cup
    \left(\bigcup_{\ell=2}^\ga T_\ell\right)
    =\left(\bigcup_{j=2}^{c_1} \gY_1(\Tone_j)\right) \cup
    \left(\bigcup_{\ell=2}^{\ga} \bigcup_{j=1}^{c_1}\gY_\ell(\Tone_j)\right).
  \end{equation}
  For each $1\leq \ell\leq \ga-1$, $\gY_\ell(\Tone_{c_1})\cup
  \gY_{\ell+1}(\Tone_1)=\Ttwo_\ell$. Furthermore, other
  $\gY_1(\Tone_j)$ and $\gY_\ell(\Tone_j)$ in the right-hand side of (\ref{eq:A-decomp}) are $T$-separate. By Remark~\ref{rmk:sT-map}, it is easy to see that
  $L_k(T_\i)\setminus L_{k+1}(T_\i)\in \sT^*$.
\eproof

\medskip

From this fact, we have the following lemma.
\begin{lemma}\label{lem:gds-TD-2}
  For any $i\in\gS_T$, $\Ttwo_i$ and $\Dtwo_i$ have same dust-like decomposition w.r.t. $\sT$
  and $\sD$.
\end{lemma}
\Proof For each touching letter $i$, we have following disjoint unions.
\begin{align*}
  &\Ttwo_i=\big(R_0(T_i)\setminus R_q(T_i)\big) \cup
  \big(L_0(T_{i+1})\setminus L_{p}(T_{i+1})\big) \cup
  \Tthree_i, \quad \mbox{and} \\
  &\Dtwo_i=\big(R_0(D_i)\setminus R_q(D_i)\big) \cup
  \big(L_0(D_{i+1})\setminus L_{p}(D_{i+1})\big) \cup \Dthree_i.
\end{align*}
The lemma follows from Lemma~\ref{lem:LuLv-SC}. \eproof


\medskip

Given $\i=i_1i_2\cdots i_m$, $\j=j_1j_2\ldots j_m\in \gS_n^*$ with the same length. We denote by $\i<\j$ if there exists $1\leq k\leq m$ such that $i_k<j_k$ and $i_t=j_t$ for $1\leq t<k$. We denote by $\i\leq \j$ if $\i<\j$ or $\i=\j$.

Given $\i,\j\in \gS_n^*$ with $\i<\j$. We say that $(\i,\j)$ is a \emph{joint pair} if $\k\leq \i$ for every $\k\in \gS_n^*$ with $\k<\j$.

\begin{lemma}\label{lem:int-Tstar}
  Let $\i,\j\in \gS_n^*$ with $\gY_\i(0)< \gY_\j(1)$. Suppose that $[\gY_\i(0),\gY_\j(1)] \cap T$ is $T$-separate and $\max\{|\i|,|\j|\}\leq \min\{p,q\}$. Then $[\gY_\i(0),\gY_\j(1)] \cap T \in \sT^*$.
\end{lemma}
\Proof
  Given $\k\in\gS_n^*$, we define the middle part of $T_\k$ to be
    $M(T_\k)=\bigcup_{j=\ga+1}^{n-\gb} T_{\k j}.$
  It is clear that $M(T_\k)=T_\k\setminus \{L_0(T_\k) \cup R_0(T_\k)\}$. We remark that $M(T_\k)=\emptyset$ for all $\k$ if $\ga=n-\gb$. In case that $\ga<n-\gb$, it is clear that $M(T_\k)\in \sT^*$ for all $\k$.

  Without loss of generality, we may assume that $|\i|\leq |\j|$. Define $k=|\j|-|\i|$. Then $[\gY_{\i [1]^k}(0), \gY_{\j}(1)] \cap T$ is $T$-separate since $\gY_{\i[1]^k}(0)=\gY_\i(0)$. Thus, noticing that the lemma holds in case that $\i=\j$, we assume that $\i<\j$ in the sequel of the proof.

  Define $m=|\i|$. Let $\gS(\i,\j)=\{\k\in\gS_n^m:\, \i\leq \k\leq \j\}$. Then
  \begin{equation}\label{eq:gYi0gYj1-decomp}
    [\gY_\i(0),\gY_\j(1)] \cap T =\bigcup_{\k\in\gS(\i,\j)} T_\k.
  \end{equation}

  Now we arbitrary pick a joint pair $(\u,\v)$ with $\u,\v\in \gS(\i,\j)$. Notice that
  \begin{equation*}
    T_\u=L_0(T_\u) \cup R_0(T_\u) \cup M(T_\u), \quad
    T_\v=L_0(T_\v) \cup R_0(T_\v) \cup M(T_\v).
  \end{equation*}

  In case that $R_0(T_\u)$ is $T$-separate, we have $R_0(T_\u)\in\sTOne$. Also, in this case, we must have $L_0(T_\v)$ is $T$-separate so that $L_0(T_\v)\in\sTOne$. It follows that there exists $A(\u,\v)\in \sT^*$ such that
  \begin{equation}\label{eq:TuTv-decomp}
    T_\u \cup T_\v = L_0(T_\u) \cup R_0(T_\v) \cup A(\u,\v),
  \end{equation}
  where the union is disjoint.

  In case that $R_0(T_\u)$ is not $T$-separate, we define $s$ to be the maximal nonnegative integer which satisfies $\u=\u^\prime[n]^s$ for some $\u^\prime\in\gS_n^*$. Let $\u^\prime=u_1u_2\cdots u_{m-s}$. From Lemma~\ref{lem:gY-i-Tonej-sep}, we have $u_{m-s}\in \gS_T$ so that $\v=\v^\prime [1]^s$ where $\v^\prime=u_1\cdots u_{m-s-1}(u_{m-s}+1)$. It is clear that $s<\min\{p,q\}$ since $|\u| =|\i|\leq \min\{p,q\}$. Notice that
  \begin{equation*}
    R_0(T_\u)=R_{q-s}(T_\u) \cup \Big(R_0(T_\u) \setminus R_{q-s}(T_\u)\Big), \quad
    L_0(T_\v)=L_{p-s}(T_\v) \cup \Big(L_0(T_\v) \setminus L_{p-s}(T_\v)\Big),
  \end{equation*}
  where the unions are disjoint and $R_0(T_\u) \setminus R_{q-s}(T_\u), L_0(T_\v) \setminus L_{p-s}(T_\v) \in \sT^*$ by Lemma~\ref{lem:LuLv-SC}. Since
  \begin{align*}
    R_{q-s}(T_\u)\cup L_{p-s}(T_\v) &= R_{q-s}(T_{\u^\prime [n]^s}) \cup L_{p-s} (T_{\v^\prime [1]^s})
    = R_q (T_{\u^\prime}) \cup L_p (T_{\v^\prime}) \\
    &= \gY_{u_1\cdots u_{m-s-1}} \Big(R_q( T_{u_{m-s}} ) \cup L_p( T_{u_{m-s}+1} ) \Big)
    = \gY_{u_1\cdots u_{m-s-1}} (\Tthree_{k_{m-s}}) \in \sTThree,
  \end{align*}
  we know that in this case, there also exists $A(\u,\v)\in \sT^*$ such that (\ref{eq:TuTv-decomp}) holds while the union is disjoint.

   Notice that $L_0(T_\i)$ and $ R_0(T_\j)$ are $T$-separate so that they are all in $\sT^*$. Using (\ref{eq:gYi0gYj1-decomp}) and (\ref{eq:TuTv-decomp}), we can see that the lemma holds.
\eproof

\begin{corol}\label{corol:Rq-R3q-Lk}
  Given $i=1,2,\ldots,n$, $k\in \bN$ and $\j\in \gS_n^*$ with $k+|\j|<\min\{p,q\}$. Assume that the last letter of $\j$ does not belong to $\{1\}\cup (\gS_T+1)$. Then $R_q(T_i)\setminus \Big(R_{3q}(T_i) \cup L_k(T_{i[n]^{2q}\j})\Big) \in \sT^*$.
\end{corol}
\Proof
  Let $\j=j_1j_2\cdots j_m$. Then $j_m>1$. Define $\u=j_1\cdots j_{m-1}(j_m-1)$. It is easy to check that
  \begin{align*}
    R_q(T_i)\setminus \Big(R_{3q}(T_i) \cup L_k(T_{i[n]^{2q}\j})\Big) = \Big([a_1, b_1] \cup [a_2, b_2] \Big) \cap T ,
  \end{align*}
  where $a_1=\gY_{i[n]^q(n-\gb+1)}(0)$, $b_1=\gY_{i[n]^{2q}\u}(1)$, $a_2=\gY_{i[n]^{2q}\j [1]^k(\ga+1)}(0)$, $b_2=\gY_{i[n]^{3q}(n-\gb)}(1)$. Thus it suffices to show that $[a_1,b_1]\cap T\in\sT^*$ and $[a_2,b_2]\cap T\in \sT^*$. Notice that
  \begin{align*}
    [a_1,b_1]\cap T = \Big(R_q(T_i)\setminus R_{2q-1}(T_i)\Big) \cup \Big([a_1^\prime, b_1]\cap T\Big),
  \end{align*}
  where $a_1^\prime=\gY_{i[n]^{2q-1}(n-\gb+1)}(0)$ and the union is disjoint. Using Lemma~\ref{lem:int-Tstar}, we have
  \begin{align*}
    [a_1^\prime, b_1]\cap T
    =\gY_{i[n]^{2q-1}}\Big( [\gY_{(n-\gb+1)}(0), \gY_{n\u}(1)] \cap T\Big) \in \sT^*
  \end{align*}
  so that $[a_1,b_1] \cap T \in \sT^*$.

  Let $s$ be the maximal nonnegative integer such that $\j=[n]^s \j^\prime$ for some $\j^\prime\in \gS_n^*$. Then $s\leq |\j|-1<q-1$. Notice that
  \begin{align*}
    [a_2,b_2]\cap T = \Big(R_{2q+s+1}(T_i)\setminus R_{3q}(T_i)\Big) \cup \Big([a_2, b_2^\prime]\cap T\Big),
  \end{align*}
  where $b_2^\prime=\gY_{i[n]^{2q+s+1}(n-\gb)}(1)$ and the union is disjoint. From $|\j^\prime| +k+1\leq \min\{p,q\}$ and Lemma~\ref{lem:int-Tstar}, we have
  \begin{align*}
    [a_2, b_2^\prime]\cap T
    =\gY_{i[n]^{2q+s}}\Big( [\gY_{\j^\prime [1]^{k}(\ga+1)}(0), \gY_{n(n-\gb)}(1)] \cap T\Big) \in \sT^*
  \end{align*}
  so that $[a_2,b_2] \cap T \in \sT^*$.
\eproof

\medskip

 The following lemma is useful in the proof of
Lemma~\ref{lem:gds-TD-34}.

\begin{lemma}\label{lem:same-set-map-TD-4}
  For any left substitutable touching letter $i$, we have
  \begin{align}
    \gY_{i[n]^{2q}}\circ \gY_{i}^{-1}(\Tfour_i)=R_{3q}(T_i) \cup
  L_{2p+k_i}(T_{i+1}), \label{eq:lem-TD-4-same-1}\\
    \gF_{i[n]^{2q}}\circ\gF_{i}^{-1}(\Dfour_i)=
  L_{k_i^\prime}(D_{i[n]^{2q}\j_i}) \cup R_{3q}(D_i).
   \label{eq:lem-TD-4-same-2}
  \end{align}
\end{lemma}
\Proof By definition of $\Tfour_i$, in order to
prove~(\ref{eq:lem-TD-4-same-1}), it suffices to
  show that
  \begin{equation*}
    \gY_{i[n]^{2q}}\circ \gY_{i}^{-1}(R_q(T_i))=R_{3q}(T_i) \quad
    \mbox{and} \quad \gY_{i[n]^{2q}}\circ
    \gY_{i}^{-1}(L_{k_i}(T_{i+1}))=L_{2p+k_i}(T_{i+1}).
  \end{equation*}

  It is clear that $$\diam\gY_{i[n]^{2q}}\circ \gY_{i}^{-1}(R_q(T_i))=\diam
  R_{3q}(T_i)\quad \mbox{and}$$
  $$\diam\gY_{i[n]^{2q}}\circ
    \gY_{i}^{-1}(L_{k_i}(T_{i+1}))=\diam L_{2p+k_i}(T_{i+1}).$$

  Notice that the maximum value of $\gY_{i[n]^{2q}}\circ
  \gY_{i}^{-1}(R_q(T_i))$ is
  \begin{equation*}
    \gY_{i[n]^{2q}}\circ \gY_{i}^{-1}(\gY_{i}(1))=\gY_{i
    [n]^{2q}}(1)=\gY_{i}(1),
  \end{equation*}
  which equals the maximum value of $R_{3q}(T_i)$. It follows
  that $\gY_{i[n]^{2q}}\circ \gY_{i}^{-1}(R_q(T_i))=R_{3q}(T_i)$.

  Since $i$ is a touching letter, the minimum value of
  $T_{i+1}$ equals $\gY_{i}(1)$. Thus the minimum
  value of $\gY_{i[n]^{2q}}\circ
  \gY_{i}^{-1}(L_{k_i}(T_{i+1}))$ is also $\gY_{i}(1)$, which
  is equals the minimum value of $L_{2p+k_i}(T_{i+1})$. If follows
  that $\gY_{i[n]^{2q}}\circ
  \gY_{i}^{-1}(L_{k_i}(T_{i+1}))=L_{2p+k_i}(T_{i+1}).$

  In order to prove (\ref{eq:lem-TD-4-same-2}), it suffices to show
  that
  \begin{equation*}
    \gF_{i[n]^{2q}}\circ \gF_{i}^{-1}(L_{k_i^\prime}(D_{i\j_i}))
    =L_{k_i^\prime}(D_{i[n]^{2q}\j_i}) \quad
    \mbox{and} \quad \gF_{i[n]^{2q}}\circ
    \gF_{i}^{-1}(R_q(D_i))=R_{3q}(D_i).
  \end{equation*}
  Similarly as above, we can easily see that $\gF_{i[n]^{2q}}\circ
    \gF_{i}^{-1}(R_q(D_i))=R_{3q}(D_i)$.
  Using diameter and noticing that the minimum value of
  $\gF_{i[n]^{2q}}\circ \gF_{i}^{-1}(L_{k_i^\prime}(D_{i\j_i}))$ is
  \begin{equation*}
    \gF_{i[n]^{2q}}\circ \gF_{i}^{-1}(
    \gF_{i\j_i}(0))=\gF_{i[n]^{2q}\j_i}(0),
  \end{equation*}
  which equals the minimum value of
  $L_{k_i^\prime}(D_{i[n]^{2q}\j_i})$, we know that
  $\gF_{i[n]^{2q}}\circ \gF_{i}^{-1}(L_{k_i^\prime}(D_{i\j_i}))
  =L_{k_i^\prime}(D_{i[n]^{2q}\j_i})$.
\eproof

\medskip

The following lemma is natural.

\begin{lemma}\label{lem:gds-LuTiLvDj}
  Given $\i,\j\in\gS_n^*$. If there exist $u,v\in\bZ^+$ such that
  $L_u(T_\i)$ is $T$-separate and $\diam L_u(T_\i)=\diam L_v(T_\j)$,
  then $L_u(T_\i)$ and $L_v(D_\j)$ have same dust-like decomposition
  w.r.t. $\sT$ and $\sD$.
\end{lemma}
\Proof It is clear that
\begin{align*}
  L_u(T_\i)=\bigcup_{k=1}^\ga T_{\i [1]^u k} = \gY_{\i [1]^u}
  (\Tone_1), \qquad
  L_v(D_\j)=\bigcup_{k=1}^\ga D_{\j [1]^v k} = \gF_{\j [1]^v}
  (\Done_1).
\end{align*}
By $\diam L_u(T_\i)=\diam L_v(T_\j)$, we have
$\gr_1^u\gr_\i=\gr_1^v\gr_\j$ so that the contraction ratios of
$\gF_{\i [1]^u}$ and $\gF_{\j [1]^v}$ are same. Thus the lemma follows from that $\Tone_1$ and $\Done_1$ have same dust-like decomposition w.r.t. $\sT$ and $\sD$. \eproof

\medskip

Based on the above lemmas, now we can prove the following crucial lemma.

\begin{lemma}\label{lem:gds-TD-34}
  For any $i\in\gS_T$, the pairs $\Tthree_i$ and $\Dthree_i$, $\Tfour_i$ and $\Dfour_i$
  have same dust-like decomposition w.r.t. $\sT$ and $\sD$.
\end{lemma}
\Proof Without loss of generality, we only show that the lemma holds for
every left substitutable touching letter $i$.
By Lemma~\ref{lem:same-set-map-TD-4}, we have
\begin{align}
  &\Tthree_i=A_1
  \cup \Big(L_p(T_{i+1})\setminus L_{2p+k_i}(T_{i+1})\Big)
  \cup A_2
   \cup
  \Big(\gY_{i[n]^{2q}}\circ\gY_i^{-1}(\Tfour_i)\Big),\label{eq:Tthree-decomp}\\
  &\Dthree_i=B_1
  \cup \Big(L_p(D_{i+1})\setminus L_{2p+k_i}(D_{i+1})\Big)
  \cup B_2
  \cup
  \Big(\gF_{i[n]^{2q}}\circ\gF_i^{-1}(\Dfour_i)\Big), \quad \mbox{where}\label{eq:Dthree-decomp}
\end{align}
$$A_1=R_q(T_i)\setminus
\Big(L_{k_i^\prime}(T_{i[n]^{2q}\j_i})\cup R_{3q}(T_i)\Big), \quad
A_2=L_{k_i^\prime}(T_{i[n]^{2q}\j_i}),$$
$$B_1=R_q(D_i)\setminus
\Big(L_{k_i^\prime}(D_{i[n]^{2q}\j_i})\cup R_{3q}(D_i)\Big), \quad
B_2=L_{2p+k_i}(D_{i+1}).$$

Notice that $L_{k_i^\prime}(D_{i[n]^{2q}\j_i}) \cap R_{3q}(D_i)=\emptyset$ by (\ref{eq:pq-left-restriction}). Since $D$ is dust-like, it is clear that the union in (\ref{eq:Dthree-decomp}) is disjoint. By definition, the last letter of $\j_i$ does not belong to $\{1\}\cup (\gS_T+1)$. Thus, using Lemma~\ref{lem:gY-i-Tonej-sep}, we know that $A_2=L_{k_i^\prime}(T_{i[n]^{2q}\j_i})
=\gY_{i[n]^{2q}\j_i[1]^{k_i^\prime}}(\Tone_1)$ is $T$-separate. Hence, by $L_{k_i^\prime}(T_{i[n]^{2q}\j_i})  \cap
R_{3q}(T_i)=\emptyset$, we know that the union in (\ref{eq:Tthree-decomp}) is disjoint.

By Lemma~\ref{lem:gds-LuTiLvDj}, we know that $L_p(T_{i+1})\setminus L_{2p+k_i}(T_{i+1})$ and $L_p(D_{i+1})\setminus L_{2p+k_i}(D_{i+1})$ have same dust-like decomposition w.r.t. $\sT$ and $\sD$.
Thus in order to show that $\Tthree_i$ and $\Dthree_i$ have
same dust-like decomposition w.r.t. $\sT$ and $\sD$, it suffices to show
that $A_1$ and $B_1$,
$A_2$ and $B_2$ have same dust-like decomposition w.r.t. $\sT$ and $\sD$.

Notice that $|\j_i|+k_i^\prime<\min\{p,q\}$ and $L_{k_i^\prime}(T_{i[n]^{2q}\j_i})$ is $T$-separate. By Corollary~\ref{corol:Rq-R3q-Lk}, we have
$A_1\in \sT^*$. It is clear that
$B_1=\gp_D\circ \gp_T^{-1} (A_1)$.
By Lemma~\ref{lem:gpDgpTInv-decom}, $A_1$ and $B_1$
have same dust-like decomposition w.r.t. $\sT$ and $\sD$.

Notice that $\diam A_2=\gr_n^{2q}\diam L_{k_i^\prime}(T_{i\j_i})$
and $\diam L_{2p+k_i}(T_{i+1}) = \gr_1^{2p}\diam L_{k_i}(T_{i+1})$.
From the definition of $k_i,k_i^\prime$ and $\j_i$, we know that
$\diam L_{k_i^\prime}(T_{i\j_i})=\diam L_{k_i}(T_{i+1})$ so that
$\diam A_2=\diam L_{2p+k_i}(T_{i+1})$. Since $A_2$ is $T$-separate, by Lemma~\ref{lem:gds-LuTiLvDj}, $A_2$ and $B_2$ have
same dust-like decomposition w.r.t. $\sT$ and $\sD$. Hence, $\Tthree_i$ and $\Dthree_i$ have same dust-like decomposition w.r.t. $\sT$ and $\sD$.


Using Lemma~\ref{lem:same-set-map-TD-4} again,
\begin{align*}
  &\Tfour_i=A_1  \cup A_3  \cup A_2   \cup
  \Big(\gY_{i[n]^{2q}}\circ\gY_i^{-1}(\Tfour_i)\Big), \\
  &\Dfour_i=B_1  \cup B_3  \cup B_2^\prime  \cup
  \Big(\gF_{i[n]^{2q}}\circ\gF_i^{-1}(\Dfour_i)\Big), \quad\mbox{where}
\end{align*}
\begin{align*}
  A_3=L_{k_i}(T_{i+1})\setminus L_{2p+k_i}(T_{i+1}), \quad
  B_3=L_{k_i^\prime}(D_{i\j_i})
  \setminus L_{2p+k_i^\prime}(D_{i\j_i}),
  \quad B_2^\prime=L_{2p+k_i^\prime}(D_{i\j_i}).
\end{align*}
Similarly as above, we can see that $A_2$ and $B_2^\prime$ have same
dust-like decomposition w.r.t. $\sT$ and $\sD$.

By Lemma~\ref{lem:gds-LuTiLvDj}, $A_3$ and $L_{k_{i}}(D_{i+1})\setminus
L_{2p+k_i}(D_{i+1})$ have same
dust-like decomposition w.r.t. $\sT$ and $\sD$. Notice that
\begin{align*}
  L_{k_{i}}(D_{i+1})\setminus L_{2p+k_i}(D_{i+1})
  &=\gF_{(i+1)[1]^{k_i}} \Big(\Done_1\setminus
   \gF_{[1]^{2p}} (\Done_1) \Big) , \\
  B_3&=\gF_{i \j_i [1]^{k_i^\prime}} \Big(\Done_1\setminus
   \gF_{[1]^{2p}} (\Done_1) \Big).
\end{align*}
By the definition of $k_i,k_i^\prime$ and $\j_i$, we know that
$\gr_{i+1}\gr_1^{k_i}=\gr_{1}^{k_i^\prime}\gr_{i\j_i}$. Thus it is easy to see $A_3$
and $B_3$ have same decomposition w.r.t. $\sT$ and $\sD$. As a result,  $\Tfour_i$ and $\Dfour_i$ have same dust-like decomposition
w.r.t. $\sT$ and $\sD$. \eproof

\medskip

\noindent{\bf Proof of Theorem~\ref{thm:sufficient-cond}} \quad The
theorem follows from Lemmas~\ref{lem:gds-weak},~\ref{lem:gds-TD-01}, ~\ref{lem:gds-TD-2}
and~\ref{lem:gds-TD-34}. \eproof

\medskip

 {\noindent{\bf{Acknowlegements:}} This work was partially
finished during the period when the authors visited the Morningside
Center of Mathematics, Chinese Academy of Sciences. The authors wish
to thank Prof. Zhi-Ying Wen for his invitation. The research of Ruan is supported by NSF of Zhejiang Province of China (No. Y6110128). The research of Xi is supported by NSFC (No. 11071224) and NCET of China.}


\noindent  Department of Mathematics, Zhejiang University, Hangzhou,
310027, China, \\ \indent ruanhj@zju.edu.cn.

\medskip

\noindent  Department of Mathematics, Michigan State University,
East Lansing MI, 48824, USA, \\ \indent ywang@math.msu.edu.

\medskip

\noindent  Institute of Mathematics, Zhejiang Wanli University,
Ningbo, 315100, China,\\ \indent xilf@zwu.edu.cn.

\end{document}